\documentclass[10pt]{article}
\usepackage{amsmath,amsfonts,amssymb,amsthm}
\usepackage{geometry}
\usepackage{enumitem}
\usepackage{fancyhdr}
\usepackage{bbm}
\usepackage{hyperref}
\usepackage{comment}
\usepackage{color}
\usepackage{tikz}
\usepackage{xurl}
\hypersetup{colorlinks=false}
\usepackage[
backend=biber,
style=numeric,
sorting=nyt
]{biblatex}

\allowdisplaybreaks

\newtheorem{theorem}{Theorem}[section]
\newtheorem{corollary}[theorem]{Corollary}
\newtheorem{lemma}[theorem]{Lemma}
\newtheorem{proposition}[theorem]{Proposition}

\theoremstyle{definition}
\newtheorem{definition}[theorem]{Definition}

\newtheorem{remark}[theorem]{Remark}

\renewcommand{\k}{\mathbbm{k}}

\newcommand{\Q}{\mathbb{Q}}

\setlist[enumerate]{label=\arabic*),ref=\arabic*)}

\title{$D_n$-quasi-invariant polynomials in characteristic $p$}
\author{Frank Wang}
\date{\vspace{-5ex}}

\newcommand{\triv}{\mathrm{triv}}

\newcommand{\sgn}{\mathrm{sign}}
\newcommand{\Irrep}{\mathrm{Irrep}}
\renewcommand{\char}{\text{char\,}}
\newcommand{\C}{\mathbb{C}}

\newcommand{\Z}{\mathbb{Z}}

\newcommand{\h}{\mathfrak{h}}
\newcommand{\p}{\mathfrak{p}}
\newcommand{\xto}[1]{\xrightarrow{#1}}
\renewcommand{\H}{\mathcal{H}}
\DeclareMathOperator{\Fr}{Fr}
\newcommand{\m}{\mathfrak{m}}
\renewcommand{\p}{\mathfrak{p}}

\begin{document}
\maketitle
\begin{abstract}
    We give explicit descriptions of the $m$-harmonic polynomials for dihedral groups via hypergeometric functions. Using these, we compute the parameters $m$ for which the Hilbert series of the $D_n$-quasi-invariant polynomials in characteristic $p$ differs from the Hilbert series of $D_n$-quasi-invariant polynomials in characteristic $0$. We also give explicit descriptions of the generators of the $D_n$-quasi-invariant polynomials in characteristic $p$ when this is the case.
\end{abstract}

\section{Introduction}

Let $W$ be a finite Coxeter group and let $\h$ be its reflection representation. Let $(,)$ be a $W$-invariant form on $\h^*$. Let $\Sigma\subset W$ be the set of reflections, and for each $s\in\Sigma$, let $\alpha_s\in\h^*$ be the unique vector (up to sign) such that $s\alpha_s=-\alpha_s$ and $(\alpha_s,\alpha_s)=1$.

\begin{definition}
    Let $\k$ be a field, let $W$ be a finite Coxeter group, and let $m:\Sigma\to\Z_{\geq 0}$ be a function that is invariant on conjugacy classes. Let $\h$ be the reflection representation of $W$. A polynomial $P\in\k[\h]$ is called \textbf{$m$-quasi-invariant} if for all $s\in\Sigma$, $(1-s)P$ is divisible by $\alpha_s(x)^{2m_s+1}$. We denote by $Q_m(W,\k)$ the space of $m$-quasi-invariants.
\end{definition}

Note that $Q_m(W,\k)$ is both a ring and a module over $\k[\h]^W$, e.g.~see \cite{etingof2002lectures}. $Q_m(W,\k)$ is also preserved by the action of $W$, hence if $\char\k\nmid|W|$ we have a decomposition
\[Q_m(W,\k)=\bigoplus_{\tau\in\Irrep(W)}Q_m(W,\k)_\tau\]
of $Q_m(W,\k)$ into isotypic components with respect to the $W$-action. 

For a graded vector space $V=\bigoplus_{d\geq 0}V_d$, let
\[\H(V)=\sum_{d\geq 0}\dim(V_d)t^d\]
be its Hilbert series. We are interested in computing the Hilbert series $\H(Q_m(W,\k))$. Note that since the equations defining quasi-invariant polynomials are defined over the integers, $\H(Q_m(W,\k))$ only depends on the characteristic of $\k$.

When $\char\k=0$, it was proved in \cite{eg} that $Q_m(W,\k)$ is free over $\k[\h]^W$, and the degrees of its generators were computed in \cite{felder2003action}. Explicit generators of the quasi-invariants for dihedral groups $D_n$ were computed in \cite{feigin}, \cite{feigindihedral}. In particular, the Hilbert series of the space of quasi-invariants in characteristic $0$ is known. However, in characteristic $p$, the Hilbert series is still unknown in general, and there has been work in recent years to compute it (see \cite{Ren_2020}, \cite{wang}, \cite{wangyee}, \cite{judy}). In this paper we prove freeness of $Q_m(D_n,\k)$ over a field of characteristic $p$, along with an explicit description of its generators. In particular, we have a full description of $\H(Q_m(D_n,\k))$ and conditions for when it differs from the characteristic $0$ Hilbert series.

Quasi-invariant polynomials play a key role in the representation theory of rational Cherednik algebras (see \cite{beg}, \cite{etingof2002lectures}), where they contain Verma modules as summands. They are also closely related to $m$-harmonic polynomials (see \cite{feigin}, \cite{eg}), which are the joint kernel of the quantum integrals of the quantum Calogero-Moser problem. This latter perspective will be crucial to this paper, where a large step in the proof of the main theorem is an explicit description of the $m$-harmonic polynomials for dihedral groups in terms of hypergeometric functions.

We will use the following notation for $D_n$: $D_n$ has the identity element $1$, reflections $s_0,s_1,\dots,s_{n-1}$, and rotations $t,t^2,\dots,t^{n-1}$, with relations
\[t^n=1,\quad ts_j=s_{j-1},\quad s_jt=s_{j+1},\quad s_j^2=1.\]
When $n$ is odd, there is a single conjugacy class of reflections, thus the multiplicity function $m$ is constant, and we will denote its value on any reflection also by $m$. When $n$ is even, there are two conjugacy classes of reflections given by $\{s_0,s_2,\dots,s_{n-2}\}$ and $\{s_1,s_3,\dots,s_{n-1}\}$, and we will denote the value of $m$ on them by $m_0$ and $m_1$, respectively. We will also let
\[\mu:=m_0+m_1,\quad\delta=m_0-m_1.\]

Let $\k$ be a field with $\char\k\nmid2n$ such that $\k$ has $n$-th roots of unity. Fix a primitive $n$-th root of unity in $\k$ and call it $e^{2\pi i/n}$. Then the action of $D_n$ on $\h\cong\h^*$ is given by
\[s_j\mapsto\begin{bmatrix}
    0&e^{-2j\pi i/n}\\e^{2j\pi i/n}&0
\end{bmatrix},\quad t\mapsto\begin{bmatrix}
    e^{2\pi i/n}&0\\0&e^{-2\pi i/n}
\end{bmatrix}.\]
Writing the corresponding basis of $\h^*$ in polar coordinates as $\{re^{i\theta},re^{-i\theta}\}$, we see that the action of $s_j$ corresponds to the map $\theta\mapsto\frac{2j\pi}{n}-\theta$ and the action of $t$ corresponds to the map $\theta\mapsto\theta+\frac{2\pi}{n}$. Then we have $\alpha_j(x):=\alpha_{s_j}(x)=r(e^{i\theta}-e^{-i\theta+2j\pi i/n})$.

We denote by $\triv,\sgn$ the trivial and sign representations of $D_n$. For $1\leq k<\frac{n}{2}$, we let $\rho_k$ be the two-dimensional representation whose action is given by
\[s_j\mapsto\begin{bmatrix}
    0&e^{-2kj\pi i/n}\\e^{2kj\pi i/n}&0
\end{bmatrix},\quad t\mapsto\begin{bmatrix}
    e^{2k\pi i/n}&0\\0&e^{-2k\pi i/n}
\end{bmatrix}.\]
Let the two basis vectors of $\rho_k$ with respect to this matrix presentation be $v_k$ and $s_0v_k$, respectively.

Note that we have $\rho_1=\h^*$. When $n$ is even, we also have two additional one-dimensional representations which we denote by $\rho_{\pm}$ and $\rho_{\mp}$; the action of $s_j$ on $\rho_{\pm}$ is by $(-1)^j$ and the action of $s_j$ on $\rho_{\mp}$ is by $(-1)^{j+1}$.

For the one-dimensional representations $\tau$, the Hilbert series $\H(Q_m(D_n,\k)_\tau)$ were computed over all characteristics in \cite{feigindihedral}, \cite{judy}. So in this paper we will focus on the two-dimensional representations $\rho_k$. 

Let $\char\k=p$. The Frobenius map $\Fr:\k[D_n]\to\k[D_n]$ induces a pullback map $\Fr^*:\Irrep(\k[D_n])\to\Irrep(\k[D_n])$. Then we define a Frobenius map $\Fr:\{1,\dots,\lfloor\frac{n-1}{2}\rfloor\}\to\{1,\dots,\lfloor\frac{n-1}{2}\rfloor\}$ such that if $1\leq k,k'<\frac{n}{2}$ are such that $\Fr^*\rho_{k'}\cong\rho_k$, then $\Fr(k)=k'$. Explicitly, we have
\begin{equation}\label{eq:frob}
    \Fr(k)=\min(kp^{-1}\mod n, -kp^{-1}\mod n)
\end{equation}
where $x\mod n$ is the unique integer between $0$ and $n-1$ inclusive that is congruent to $x$ modulo $n$.

For $x\in\Q$, let $\{x\}$ denote the fractional part of $x$. Let $\k_0$ be a field of characteristic $0$ and $\k_p$ be a field of characteristic $p$, such that both $\k_0$ and $\k_p$ contain $n$-th roots of unity. We have the following theorem:

\begin{theorem}\label{thm:condition}
    If $n$ is odd, then the Hilbert series of $Q_m(D_n,\k_p)_{\rho_k}$ is different from the Hilbert series of $Q_m(D_n,\k_0)_{\rho_k}$ if and only if there is some $a$ such that
    \begin{equation}\label{eq:odddet}
        \frac{\Fr^a(k)}{n}<\left\{\frac{m}{p^a}\right\}<\frac{n-\Fr^a(k)}{n}-\frac{1}{p^a}.
    \end{equation}
    If $n$ is even, then the Hilbert series of $Q_m(D_n,\k_p)_{\rho_k}$ is different from the Hilbert series of $Q_m(D_n,\k_0)_{\rho_k}$ if and only if there is some $a$ such that either
    \begin{equation}\label{eq:evendet1}
        \frac{\Fr^a(k)}{n}<\left\{\frac{\mu}{2p^a}\right\}<\frac{n-\Fr^a(k)}{n}-\frac{1}{p^a},\quad\frac{\Fr^a(k)}{n}+\frac{1}{2p^a}<\left\{\frac{\delta}{2p^a}+\frac{1}{2}\right\}<\frac{n-\Fr^a(k)}{n}-\frac{1}{2p^a}
    \end{equation}
    or
    \begin{equation}\label{eq:evendet2}
        \frac{1}{2}-\frac{\Fr^a(k)}{n}<\left\{\frac{\mu}{2p^a}+\frac{1}{2}\right\}<\frac{1}{2}+\frac{\Fr^a(k)}{n}-\frac{1}{p^a},\quad\frac{1}{2}-\frac{\Fr^a(k)}{n}+\frac{1}{2p^a}<\left\{\frac{\delta}{2p^a}\right\}<\frac{1}{2}+\frac{\Fr^a(k)}{n}-\frac{1}{2p^a}.
    \end{equation}
\end{theorem}

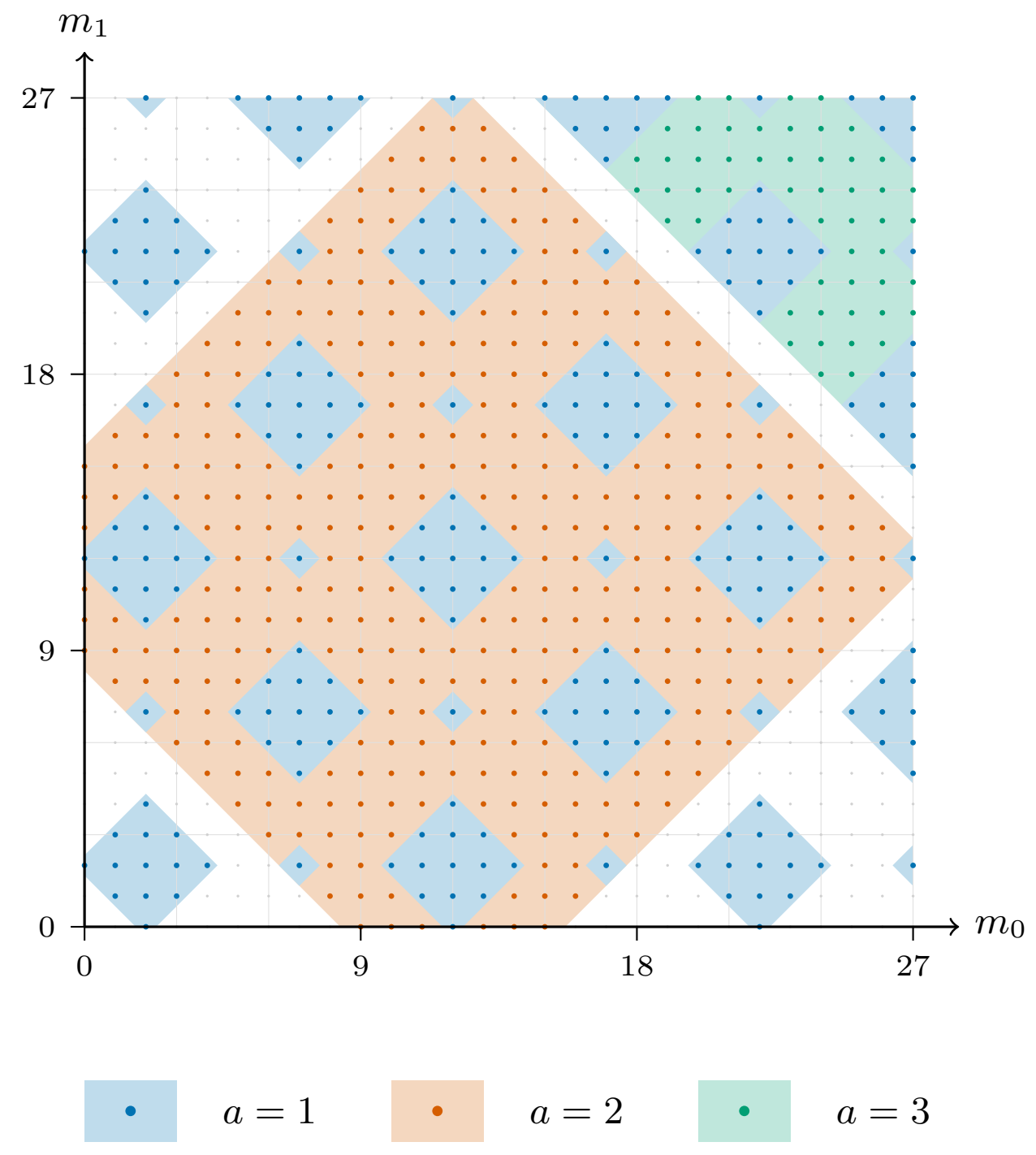
\begin{figure}
    \centering
\definecolor{qaone}{HTML}{0072B2}
\definecolor{qatwo}{HTML}{D55E00}
\definecolor{qathree}{HTML}{009E73}
\definecolor{qafour}{HTML}{CC79A7}
\begin{tikzpicture}[x=0.25cm,y=0.25cm,font=\small]
\begin{scope}[xshift=0cm]
\fill[white] (0,0) rectangle (27,27);
\fill[qathree!25] (27.000000,14.666667) -- (27.000000,27.000000) -- (14.666667,27.000000) -- cycle;
\fill[qatwo!25] (8.333333,0.000000) -- (15.666667,0.000000) -- (27.000000,11.333333) -- (27.000000,12.666667) -- (12.666667,27.000000) -- (11.333333,27.000000) -- (0.000000,15.666667) -- (0.000000,8.333333) -- cycle;
\fill[qaone!25] (1.666667,0.000000) -- (2.333333,0.000000) -- (4.333333,2.000000) -- (2.000000,4.333333) -- (0.000000,2.333333) -- (0.000000,1.666667) -- cycle;
\fill[qaone!25] (4.333333,12.000000) -- (2.000000,14.333333) -- (0.000000,12.333333) -- (0.000000,11.666667) -- (2.000000,9.666667) -- cycle;
\fill[qaone!25] (9.333333,7.000000) -- (7.000000,9.333333) -- (4.666667,7.000000) -- (7.000000,4.666667) -- cycle;
\fill[qaone!25] (11.666667,0.000000) -- (12.333333,0.000000) -- (14.333333,2.000000) -- (12.000000,4.333333) -- (9.666667,2.000000) -- cycle;
\fill[qaone!25] (4.333333,22.000000) -- (2.000000,24.333333) -- (0.000000,22.333333) -- (0.000000,21.666667) -- (2.000000,19.666667) -- cycle;
\fill[qaone!25] (9.333333,17.000000) -- (7.000000,19.333333) -- (4.666667,17.000000) -- (7.000000,14.666667) -- cycle;
\fill[qaone!25] (14.333333,12.000000) -- (12.000000,14.333333) -- (9.666667,12.000000) -- (12.000000,9.666667) -- cycle;
\fill[qaone!25] (19.333333,7.000000) -- (17.000000,9.333333) -- (14.666667,7.000000) -- (17.000000,4.666667) -- cycle;
\fill[qaone!25] (21.666667,0.000000) -- (22.333333,0.000000) -- (24.333333,2.000000) -- (22.000000,4.333333) -- (19.666667,2.000000) -- cycle;
\fill[qaone!25] (9.333333,27.000000) -- (4.666667,27.000000) -- (7.000000,24.666667) -- cycle;
\fill[qaone!25] (14.333333,22.000000) -- (12.000000,24.333333) -- (9.666667,22.000000) -- (12.000000,19.666667) -- cycle;
\fill[qaone!25] (19.333333,17.000000) -- (17.000000,19.333333) -- (14.666667,17.000000) -- (17.000000,14.666667) -- cycle;
\fill[qaone!25] (24.333333,12.000000) -- (22.000000,14.333333) -- (19.666667,12.000000) -- (22.000000,9.666667) -- cycle;
\fill[qaone!25] (27.000000,4.666667) -- (27.000000,9.333333) -- (24.666667,7.000000) -- cycle;
\fill[qaone!25] (19.333333,27.000000) -- (14.666667,27.000000) -- (17.000000,24.666667) -- cycle;
\fill[qaone!25] (24.333333,22.000000) -- (22.000000,24.333333) -- (19.666667,22.000000) -- (22.000000,19.666667) -- cycle;
\fill[qaone!25] (27.000000,14.666667) -- (27.000000,19.333333) -- (24.666667,17.000000) -- cycle;
\fill[qaone!25] (27.000000,24.666667) -- (27.000000,27.000000) -- (24.666667,27.000000) -- cycle;
\fill[qaone!25] (2.666667,7.000000) -- (2.000000,7.666667) -- (1.333333,7.000000) -- (2.000000,6.333333) -- cycle;
\fill[qaone!25] (7.666667,2.000000) -- (7.000000,2.666667) -- (6.333333,2.000000) -- (7.000000,1.333333) -- cycle;
\fill[qaone!25] (2.666667,17.000000) -- (2.000000,17.666667) -- (1.333333,17.000000) -- (2.000000,16.333333) -- cycle;
\fill[qaone!25] (7.666667,12.000000) -- (7.000000,12.666667) -- (6.333333,12.000000) -- (7.000000,11.333333) -- cycle;
\fill[qaone!25] (12.666667,7.000000) -- (12.000000,7.666667) -- (11.333333,7.000000) -- (12.000000,6.333333) -- cycle;
\fill[qaone!25] (17.666667,2.000000) -- (17.000000,2.666667) -- (16.333333,2.000000) -- (17.000000,1.333333) -- cycle;
\fill[qaone!25] (2.666667,27.000000) -- (1.333333,27.000000) -- (2.000000,26.333333) -- cycle;
\fill[qaone!25] (7.666667,22.000000) -- (7.000000,22.666667) -- (6.333333,22.000000) -- (7.000000,21.333333) -- cycle;
\fill[qaone!25] (12.666667,17.000000) -- (12.000000,17.666667) -- (11.333333,17.000000) -- (12.000000,16.333333) -- cycle;
\fill[qaone!25] (17.666667,12.000000) -- (17.000000,12.666667) -- (16.333333,12.000000) -- (17.000000,11.333333) -- cycle;
\fill[qaone!25] (22.666667,7.000000) -- (22.000000,7.666667) -- (21.333333,7.000000) -- (22.000000,6.333333) -- cycle;
\fill[qaone!25] (27.000000,1.333333) -- (27.000000,2.666667) -- (26.333333,2.000000) -- cycle;
\fill[qaone!25] (12.666667,27.000000) -- (11.333333,27.000000) -- (12.000000,26.333333) -- cycle;
\fill[qaone!25] (17.666667,22.000000) -- (17.000000,22.666667) -- (16.333333,22.000000) -- (17.000000,21.333333) -- cycle;
\fill[qaone!25] (22.666667,17.000000) -- (22.000000,17.666667) -- (21.333333,17.000000) -- (22.000000,16.333333) -- cycle;
\fill[qaone!25] (27.000000,11.333333) -- (27.000000,12.666667) -- (26.333333,12.000000) -- cycle;
\fill[qaone!25] (22.666667,27.000000) -- (21.333333,27.000000) -- (22.000000,26.333333) -- cycle;
\fill[qaone!25] (27.000000,21.333333) -- (27.000000,22.666667) -- (26.333333,22.000000) -- cycle;
\draw[step=3,black!12,very thin] (0,0) grid (27,27);
\foreach \x in {0,...,27} \foreach \y in {0,...,27} \fill[black!18] (\x,\y) circle[radius=0.32pt];
\fill[qaone] (0,2) circle[radius=0.64pt];
\fill[qatwo] (0,9) circle[radius=0.64pt];
\fill[qatwo] (0,10) circle[radius=0.64pt];
\fill[qatwo] (0,11) circle[radius=0.64pt];
\fill[qaone] (0,12) circle[radius=0.64pt];
\fill[qatwo] (0,13) circle[radius=0.64pt];
\fill[qatwo] (0,14) circle[radius=0.64pt];
\fill[qatwo] (0,15) circle[radius=0.64pt];
\fill[qaone] (0,22) circle[radius=0.64pt];
\fill[qaone] (1,1) circle[radius=0.64pt];
\fill[qaone] (1,2) circle[radius=0.64pt];
\fill[qaone] (1,3) circle[radius=0.64pt];
\fill[qatwo] (1,8) circle[radius=0.64pt];
\fill[qatwo] (1,9) circle[radius=0.64pt];
\fill[qatwo] (1,10) circle[radius=0.64pt];
\fill[qaone] (1,11) circle[radius=0.64pt];
\fill[qaone] (1,12) circle[radius=0.64pt];
\fill[qaone] (1,13) circle[radius=0.64pt];
\fill[qatwo] (1,14) circle[radius=0.64pt];
\fill[qatwo] (1,15) circle[radius=0.64pt];
\fill[qatwo] (1,16) circle[radius=0.64pt];
\fill[qaone] (1,21) circle[radius=0.64pt];
\fill[qaone] (1,22) circle[radius=0.64pt];
\fill[qaone] (1,23) circle[radius=0.64pt];
\fill[qaone] (2,0) circle[radius=0.64pt];
\fill[qaone] (2,1) circle[radius=0.64pt];
\fill[qaone] (2,2) circle[radius=0.64pt];
\fill[qaone] (2,3) circle[radius=0.64pt];
\fill[qaone] (2,4) circle[radius=0.64pt];
\fill[qaone] (2,7) circle[radius=0.64pt];
\fill[qatwo] (2,8) circle[radius=0.64pt];
\fill[qatwo] (2,9) circle[radius=0.64pt];
\fill[qaone] (2,10) circle[radius=0.64pt];
\fill[qaone] (2,11) circle[radius=0.64pt];
\fill[qaone] (2,12) circle[radius=0.64pt];
\fill[qaone] (2,13) circle[radius=0.64pt];
\fill[qaone] (2,14) circle[radius=0.64pt];
\fill[qatwo] (2,15) circle[radius=0.64pt];
\fill[qatwo] (2,16) circle[radius=0.64pt];
\fill[qaone] (2,17) circle[radius=0.64pt];
\fill[qaone] (2,20) circle[radius=0.64pt];
\fill[qaone] (2,21) circle[radius=0.64pt];
\fill[qaone] (2,22) circle[radius=0.64pt];
\fill[qaone] (2,23) circle[radius=0.64pt];
\fill[qaone] (2,24) circle[radius=0.64pt];
\fill[qaone] (2,27) circle[radius=0.64pt];
\fill[qaone] (3,1) circle[radius=0.64pt];
\fill[qaone] (3,2) circle[radius=0.64pt];
\fill[qaone] (3,3) circle[radius=0.64pt];
\fill[qatwo] (3,6) circle[radius=0.64pt];
\fill[qatwo] (3,7) circle[radius=0.64pt];
\fill[qatwo] (3,8) circle[radius=0.64pt];
\fill[qatwo] (3,9) circle[radius=0.64pt];
\fill[qatwo] (3,10) circle[radius=0.64pt];
\fill[qaone] (3,11) circle[radius=0.64pt];
\fill[qaone] (3,12) circle[radius=0.64pt];
\fill[qaone] (3,13) circle[radius=0.64pt];
\fill[qatwo] (3,14) circle[radius=0.64pt];
\fill[qatwo] (3,15) circle[radius=0.64pt];
\fill[qatwo] (3,16) circle[radius=0.64pt];
\fill[qatwo] (3,17) circle[radius=0.64pt];
\fill[qatwo] (3,18) circle[radius=0.64pt];
\fill[qaone] (3,21) circle[radius=0.64pt];
\fill[qaone] (3,22) circle[radius=0.64pt];
\fill[qaone] (3,23) circle[radius=0.64pt];
\fill[qaone] (4,2) circle[radius=0.64pt];
\fill[qatwo] (4,5) circle[radius=0.64pt];
\fill[qatwo] (4,6) circle[radius=0.64pt];
\fill[qatwo] (4,7) circle[radius=0.64pt];
\fill[qatwo] (4,8) circle[radius=0.64pt];
\fill[qatwo] (4,9) circle[radius=0.64pt];
\fill[qatwo] (4,10) circle[radius=0.64pt];
\fill[qatwo] (4,11) circle[radius=0.64pt];
\fill[qaone] (4,12) circle[radius=0.64pt];
\fill[qatwo] (4,13) circle[radius=0.64pt];
\fill[qatwo] (4,14) circle[radius=0.64pt];
\fill[qatwo] (4,15) circle[radius=0.64pt];
\fill[qatwo] (4,16) circle[radius=0.64pt];
\fill[qatwo] (4,17) circle[radius=0.64pt];
\fill[qatwo] (4,18) circle[radius=0.64pt];
\fill[qatwo] (4,19) circle[radius=0.64pt];
\fill[qaone] (4,22) circle[radius=0.64pt];
\fill[qatwo] (5,4) circle[radius=0.64pt];
\fill[qatwo] (5,5) circle[radius=0.64pt];
\fill[qatwo] (5,6) circle[radius=0.64pt];
\fill[qaone] (5,7) circle[radius=0.64pt];
\fill[qatwo] (5,8) circle[radius=0.64pt];
\fill[qatwo] (5,9) circle[radius=0.64pt];
\fill[qatwo] (5,10) circle[radius=0.64pt];
\fill[qatwo] (5,11) circle[radius=0.64pt];
\fill[qatwo] (5,12) circle[radius=0.64pt];
\fill[qatwo] (5,13) circle[radius=0.64pt];
\fill[qatwo] (5,14) circle[radius=0.64pt];
\fill[qatwo] (5,15) circle[radius=0.64pt];
\fill[qatwo] (5,16) circle[radius=0.64pt];
\fill[qaone] (5,17) circle[radius=0.64pt];
\fill[qatwo] (5,18) circle[radius=0.64pt];
\fill[qatwo] (5,19) circle[radius=0.64pt];
\fill[qatwo] (5,20) circle[radius=0.64pt];
\fill[qaone] (5,27) circle[radius=0.64pt];
\fill[qatwo] (6,3) circle[radius=0.64pt];
\fill[qatwo] (6,4) circle[radius=0.64pt];
\fill[qatwo] (6,5) circle[radius=0.64pt];
\fill[qaone] (6,6) circle[radius=0.64pt];
\fill[qaone] (6,7) circle[radius=0.64pt];
\fill[qaone] (6,8) circle[radius=0.64pt];
\fill[qatwo] (6,9) circle[radius=0.64pt];
\fill[qatwo] (6,10) circle[radius=0.64pt];
\fill[qatwo] (6,11) circle[radius=0.64pt];
\fill[qatwo] (6,12) circle[radius=0.64pt];
\fill[qatwo] (6,13) circle[radius=0.64pt];
\fill[qatwo] (6,14) circle[radius=0.64pt];
\fill[qatwo] (6,15) circle[radius=0.64pt];
\fill[qaone] (6,16) circle[radius=0.64pt];
\fill[qaone] (6,17) circle[radius=0.64pt];
\fill[qaone] (6,18) circle[radius=0.64pt];
\fill[qatwo] (6,19) circle[radius=0.64pt];
\fill[qatwo] (6,20) circle[radius=0.64pt];
\fill[qatwo] (6,21) circle[radius=0.64pt];
\fill[qaone] (6,26) circle[radius=0.64pt];
\fill[qaone] (6,27) circle[radius=0.64pt];
\fill[qaone] (7,2) circle[radius=0.64pt];
\fill[qatwo] (7,3) circle[radius=0.64pt];
\fill[qatwo] (7,4) circle[radius=0.64pt];
\fill[qaone] (7,5) circle[radius=0.64pt];
\fill[qaone] (7,6) circle[radius=0.64pt];
\fill[qaone] (7,7) circle[radius=0.64pt];
\fill[qaone] (7,8) circle[radius=0.64pt];
\fill[qaone] (7,9) circle[radius=0.64pt];
\fill[qatwo] (7,10) circle[radius=0.64pt];
\fill[qatwo] (7,11) circle[radius=0.64pt];
\fill[qaone] (7,12) circle[radius=0.64pt];
\fill[qatwo] (7,13) circle[radius=0.64pt];
\fill[qatwo] (7,14) circle[radius=0.64pt];
\fill[qaone] (7,15) circle[radius=0.64pt];
\fill[qaone] (7,16) circle[radius=0.64pt];
\fill[qaone] (7,17) circle[radius=0.64pt];
\fill[qaone] (7,18) circle[radius=0.64pt];
\fill[qaone] (7,19) circle[radius=0.64pt];
\fill[qatwo] (7,20) circle[radius=0.64pt];
\fill[qatwo] (7,21) circle[radius=0.64pt];
\fill[qaone] (7,22) circle[radius=0.64pt];
\fill[qaone] (7,25) circle[radius=0.64pt];
\fill[qaone] (7,26) circle[radius=0.64pt];
\fill[qaone] (7,27) circle[radius=0.64pt];
\fill[qatwo] (8,1) circle[radius=0.64pt];
\fill[qatwo] (8,2) circle[radius=0.64pt];
\fill[qatwo] (8,3) circle[radius=0.64pt];
\fill[qatwo] (8,4) circle[radius=0.64pt];
\fill[qatwo] (8,5) circle[radius=0.64pt];
\fill[qaone] (8,6) circle[radius=0.64pt];
\fill[qaone] (8,7) circle[radius=0.64pt];
\fill[qaone] (8,8) circle[radius=0.64pt];
\fill[qatwo] (8,9) circle[radius=0.64pt];
\fill[qatwo] (8,10) circle[radius=0.64pt];
\fill[qatwo] (8,11) circle[radius=0.64pt];
\fill[qatwo] (8,12) circle[radius=0.64pt];
\fill[qatwo] (8,13) circle[radius=0.64pt];
\fill[qatwo] (8,14) circle[radius=0.64pt];
\fill[qatwo] (8,15) circle[radius=0.64pt];
\fill[qaone] (8,16) circle[radius=0.64pt];
\fill[qaone] (8,17) circle[radius=0.64pt];
\fill[qaone] (8,18) circle[radius=0.64pt];
\fill[qatwo] (8,19) circle[radius=0.64pt];
\fill[qatwo] (8,20) circle[radius=0.64pt];
\fill[qatwo] (8,21) circle[radius=0.64pt];
\fill[qatwo] (8,22) circle[radius=0.64pt];
\fill[qatwo] (8,23) circle[radius=0.64pt];
\fill[qaone] (8,26) circle[radius=0.64pt];
\fill[qaone] (8,27) circle[radius=0.64pt];
\fill[qatwo] (9,0) circle[radius=0.64pt];
\fill[qatwo] (9,1) circle[radius=0.64pt];
\fill[qatwo] (9,2) circle[radius=0.64pt];
\fill[qatwo] (9,3) circle[radius=0.64pt];
\fill[qatwo] (9,4) circle[radius=0.64pt];
\fill[qatwo] (9,5) circle[radius=0.64pt];
\fill[qatwo] (9,6) circle[radius=0.64pt];
\fill[qaone] (9,7) circle[radius=0.64pt];
\fill[qatwo] (9,8) circle[radius=0.64pt];
\fill[qatwo] (9,9) circle[radius=0.64pt];
\fill[qatwo] (9,10) circle[radius=0.64pt];
\fill[qatwo] (9,11) circle[radius=0.64pt];
\fill[qatwo] (9,12) circle[radius=0.64pt];
\fill[qatwo] (9,13) circle[radius=0.64pt];
\fill[qatwo] (9,14) circle[radius=0.64pt];
\fill[qatwo] (9,15) circle[radius=0.64pt];
\fill[qatwo] (9,16) circle[radius=0.64pt];
\fill[qaone] (9,17) circle[radius=0.64pt];
\fill[qatwo] (9,18) circle[radius=0.64pt];
\fill[qatwo] (9,19) circle[radius=0.64pt];
\fill[qatwo] (9,20) circle[radius=0.64pt];
\fill[qatwo] (9,21) circle[radius=0.64pt];
\fill[qatwo] (9,22) circle[radius=0.64pt];
\fill[qatwo] (9,23) circle[radius=0.64pt];
\fill[qatwo] (9,24) circle[radius=0.64pt];
\fill[qaone] (9,27) circle[radius=0.64pt];
\fill[qatwo] (10,0) circle[radius=0.64pt];
\fill[qatwo] (10,1) circle[radius=0.64pt];
\fill[qaone] (10,2) circle[radius=0.64pt];
\fill[qatwo] (10,3) circle[radius=0.64pt];
\fill[qatwo] (10,4) circle[radius=0.64pt];
\fill[qatwo] (10,5) circle[radius=0.64pt];
\fill[qatwo] (10,6) circle[radius=0.64pt];
\fill[qatwo] (10,7) circle[radius=0.64pt];
\fill[qatwo] (10,8) circle[radius=0.64pt];
\fill[qatwo] (10,9) circle[radius=0.64pt];
\fill[qatwo] (10,10) circle[radius=0.64pt];
\fill[qatwo] (10,11) circle[radius=0.64pt];
\fill[qaone] (10,12) circle[radius=0.64pt];
\fill[qatwo] (10,13) circle[radius=0.64pt];
\fill[qatwo] (10,14) circle[radius=0.64pt];
\fill[qatwo] (10,15) circle[radius=0.64pt];
\fill[qatwo] (10,16) circle[radius=0.64pt];
\fill[qatwo] (10,17) circle[radius=0.64pt];
\fill[qatwo] (10,18) circle[radius=0.64pt];
\fill[qatwo] (10,19) circle[radius=0.64pt];
\fill[qatwo] (10,20) circle[radius=0.64pt];
\fill[qatwo] (10,21) circle[radius=0.64pt];
\fill[qaone] (10,22) circle[radius=0.64pt];
\fill[qatwo] (10,23) circle[radius=0.64pt];
\fill[qatwo] (10,24) circle[radius=0.64pt];
\fill[qatwo] (10,25) circle[radius=0.64pt];
\fill[qatwo] (11,0) circle[radius=0.64pt];
\fill[qaone] (11,1) circle[radius=0.64pt];
\fill[qaone] (11,2) circle[radius=0.64pt];
\fill[qaone] (11,3) circle[radius=0.64pt];
\fill[qatwo] (11,4) circle[radius=0.64pt];
\fill[qatwo] (11,5) circle[radius=0.64pt];
\fill[qatwo] (11,6) circle[radius=0.64pt];
\fill[qatwo] (11,7) circle[radius=0.64pt];
\fill[qatwo] (11,8) circle[radius=0.64pt];
\fill[qatwo] (11,9) circle[radius=0.64pt];
\fill[qatwo] (11,10) circle[radius=0.64pt];
\fill[qaone] (11,11) circle[radius=0.64pt];
\fill[qaone] (11,12) circle[radius=0.64pt];
\fill[qaone] (11,13) circle[radius=0.64pt];
\fill[qatwo] (11,14) circle[radius=0.64pt];
\fill[qatwo] (11,15) circle[radius=0.64pt];
\fill[qatwo] (11,16) circle[radius=0.64pt];
\fill[qatwo] (11,17) circle[radius=0.64pt];
\fill[qatwo] (11,18) circle[radius=0.64pt];
\fill[qatwo] (11,19) circle[radius=0.64pt];
\fill[qatwo] (11,20) circle[radius=0.64pt];
\fill[qaone] (11,21) circle[radius=0.64pt];
\fill[qaone] (11,22) circle[radius=0.64pt];
\fill[qaone] (11,23) circle[radius=0.64pt];
\fill[qatwo] (11,24) circle[radius=0.64pt];
\fill[qatwo] (11,25) circle[radius=0.64pt];
\fill[qatwo] (11,26) circle[radius=0.64pt];
\fill[qaone] (12,0) circle[radius=0.64pt];
\fill[qaone] (12,1) circle[radius=0.64pt];
\fill[qaone] (12,2) circle[radius=0.64pt];
\fill[qaone] (12,3) circle[radius=0.64pt];
\fill[qaone] (12,4) circle[radius=0.64pt];
\fill[qatwo] (12,5) circle[radius=0.64pt];
\fill[qatwo] (12,6) circle[radius=0.64pt];
\fill[qaone] (12,7) circle[radius=0.64pt];
\fill[qatwo] (12,8) circle[radius=0.64pt];
\fill[qatwo] (12,9) circle[radius=0.64pt];
\fill[qaone] (12,10) circle[radius=0.64pt];
\fill[qaone] (12,11) circle[radius=0.64pt];
\fill[qaone] (12,12) circle[radius=0.64pt];
\fill[qaone] (12,13) circle[radius=0.64pt];
\fill[qaone] (12,14) circle[radius=0.64pt];
\fill[qatwo] (12,15) circle[radius=0.64pt];
\fill[qatwo] (12,16) circle[radius=0.64pt];
\fill[qaone] (12,17) circle[radius=0.64pt];
\fill[qatwo] (12,18) circle[radius=0.64pt];
\fill[qatwo] (12,19) circle[radius=0.64pt];
\fill[qaone] (12,20) circle[radius=0.64pt];
\fill[qaone] (12,21) circle[radius=0.64pt];
\fill[qaone] (12,22) circle[radius=0.64pt];
\fill[qaone] (12,23) circle[radius=0.64pt];
\fill[qaone] (12,24) circle[radius=0.64pt];
\fill[qatwo] (12,25) circle[radius=0.64pt];
\fill[qatwo] (12,26) circle[radius=0.64pt];
\fill[qaone] (12,27) circle[radius=0.64pt];
\fill[qatwo] (13,0) circle[radius=0.64pt];
\fill[qaone] (13,1) circle[radius=0.64pt];
\fill[qaone] (13,2) circle[radius=0.64pt];
\fill[qaone] (13,3) circle[radius=0.64pt];
\fill[qatwo] (13,4) circle[radius=0.64pt];
\fill[qatwo] (13,5) circle[radius=0.64pt];
\fill[qatwo] (13,6) circle[radius=0.64pt];
\fill[qatwo] (13,7) circle[radius=0.64pt];
\fill[qatwo] (13,8) circle[radius=0.64pt];
\fill[qatwo] (13,9) circle[radius=0.64pt];
\fill[qatwo] (13,10) circle[radius=0.64pt];
\fill[qaone] (13,11) circle[radius=0.64pt];
\fill[qaone] (13,12) circle[radius=0.64pt];
\fill[qaone] (13,13) circle[radius=0.64pt];
\fill[qatwo] (13,14) circle[radius=0.64pt];
\fill[qatwo] (13,15) circle[radius=0.64pt];
\fill[qatwo] (13,16) circle[radius=0.64pt];
\fill[qatwo] (13,17) circle[radius=0.64pt];
\fill[qatwo] (13,18) circle[radius=0.64pt];
\fill[qatwo] (13,19) circle[radius=0.64pt];
\fill[qatwo] (13,20) circle[radius=0.64pt];
\fill[qaone] (13,21) circle[radius=0.64pt];
\fill[qaone] (13,22) circle[radius=0.64pt];
\fill[qaone] (13,23) circle[radius=0.64pt];
\fill[qatwo] (13,24) circle[radius=0.64pt];
\fill[qatwo] (13,25) circle[radius=0.64pt];
\fill[qatwo] (13,26) circle[radius=0.64pt];
\fill[qatwo] (14,0) circle[radius=0.64pt];
\fill[qatwo] (14,1) circle[radius=0.64pt];
\fill[qaone] (14,2) circle[radius=0.64pt];
\fill[qatwo] (14,3) circle[radius=0.64pt];
\fill[qatwo] (14,4) circle[radius=0.64pt];
\fill[qatwo] (14,5) circle[radius=0.64pt];
\fill[qatwo] (14,6) circle[radius=0.64pt];
\fill[qatwo] (14,7) circle[radius=0.64pt];
\fill[qatwo] (14,8) circle[radius=0.64pt];
\fill[qatwo] (14,9) circle[radius=0.64pt];
\fill[qatwo] (14,10) circle[radius=0.64pt];
\fill[qatwo] (14,11) circle[radius=0.64pt];
\fill[qaone] (14,12) circle[radius=0.64pt];
\fill[qatwo] (14,13) circle[radius=0.64pt];
\fill[qatwo] (14,14) circle[radius=0.64pt];
\fill[qatwo] (14,15) circle[radius=0.64pt];
\fill[qatwo] (14,16) circle[radius=0.64pt];
\fill[qatwo] (14,17) circle[radius=0.64pt];
\fill[qatwo] (14,18) circle[radius=0.64pt];
\fill[qatwo] (14,19) circle[radius=0.64pt];
\fill[qatwo] (14,20) circle[radius=0.64pt];
\fill[qatwo] (14,21) circle[radius=0.64pt];
\fill[qaone] (14,22) circle[radius=0.64pt];
\fill[qatwo] (14,23) circle[radius=0.64pt];
\fill[qatwo] (14,24) circle[radius=0.64pt];
\fill[qatwo] (14,25) circle[radius=0.64pt];
\fill[qatwo] (15,0) circle[radius=0.64pt];
\fill[qatwo] (15,1) circle[radius=0.64pt];
\fill[qatwo] (15,2) circle[radius=0.64pt];
\fill[qatwo] (15,3) circle[radius=0.64pt];
\fill[qatwo] (15,4) circle[radius=0.64pt];
\fill[qatwo] (15,5) circle[radius=0.64pt];
\fill[qatwo] (15,6) circle[radius=0.64pt];
\fill[qaone] (15,7) circle[radius=0.64pt];
\fill[qatwo] (15,8) circle[radius=0.64pt];
\fill[qatwo] (15,9) circle[radius=0.64pt];
\fill[qatwo] (15,10) circle[radius=0.64pt];
\fill[qatwo] (15,11) circle[radius=0.64pt];
\fill[qatwo] (15,12) circle[radius=0.64pt];
\fill[qatwo] (15,13) circle[radius=0.64pt];
\fill[qatwo] (15,14) circle[radius=0.64pt];
\fill[qatwo] (15,15) circle[radius=0.64pt];
\fill[qatwo] (15,16) circle[radius=0.64pt];
\fill[qaone] (15,17) circle[radius=0.64pt];
\fill[qatwo] (15,18) circle[radius=0.64pt];
\fill[qatwo] (15,19) circle[radius=0.64pt];
\fill[qatwo] (15,20) circle[radius=0.64pt];
\fill[qatwo] (15,21) circle[radius=0.64pt];
\fill[qatwo] (15,22) circle[radius=0.64pt];
\fill[qatwo] (15,23) circle[radius=0.64pt];
\fill[qatwo] (15,24) circle[radius=0.64pt];
\fill[qaone] (15,27) circle[radius=0.64pt];
\fill[qatwo] (16,1) circle[radius=0.64pt];
\fill[qatwo] (16,2) circle[radius=0.64pt];
\fill[qatwo] (16,3) circle[radius=0.64pt];
\fill[qatwo] (16,4) circle[radius=0.64pt];
\fill[qatwo] (16,5) circle[radius=0.64pt];
\fill[qaone] (16,6) circle[radius=0.64pt];
\fill[qaone] (16,7) circle[radius=0.64pt];
\fill[qaone] (16,8) circle[radius=0.64pt];
\fill[qatwo] (16,9) circle[radius=0.64pt];
\fill[qatwo] (16,10) circle[radius=0.64pt];
\fill[qatwo] (16,11) circle[radius=0.64pt];
\fill[qatwo] (16,12) circle[radius=0.64pt];
\fill[qatwo] (16,13) circle[radius=0.64pt];
\fill[qatwo] (16,14) circle[radius=0.64pt];
\fill[qatwo] (16,15) circle[radius=0.64pt];
\fill[qaone] (16,16) circle[radius=0.64pt];
\fill[qaone] (16,17) circle[radius=0.64pt];
\fill[qaone] (16,18) circle[radius=0.64pt];
\fill[qatwo] (16,19) circle[radius=0.64pt];
\fill[qatwo] (16,20) circle[radius=0.64pt];
\fill[qatwo] (16,21) circle[radius=0.64pt];
\fill[qatwo] (16,22) circle[radius=0.64pt];
\fill[qatwo] (16,23) circle[radius=0.64pt];
\fill[qaone] (16,26) circle[radius=0.64pt];
\fill[qaone] (16,27) circle[radius=0.64pt];
\fill[qaone] (17,2) circle[radius=0.64pt];
\fill[qatwo] (17,3) circle[radius=0.64pt];
\fill[qatwo] (17,4) circle[radius=0.64pt];
\fill[qaone] (17,5) circle[radius=0.64pt];
\fill[qaone] (17,6) circle[radius=0.64pt];
\fill[qaone] (17,7) circle[radius=0.64pt];
\fill[qaone] (17,8) circle[radius=0.64pt];
\fill[qaone] (17,9) circle[radius=0.64pt];
\fill[qatwo] (17,10) circle[radius=0.64pt];
\fill[qatwo] (17,11) circle[radius=0.64pt];
\fill[qaone] (17,12) circle[radius=0.64pt];
\fill[qatwo] (17,13) circle[radius=0.64pt];
\fill[qatwo] (17,14) circle[radius=0.64pt];
\fill[qaone] (17,15) circle[radius=0.64pt];
\fill[qaone] (17,16) circle[radius=0.64pt];
\fill[qaone] (17,17) circle[radius=0.64pt];
\fill[qaone] (17,18) circle[radius=0.64pt];
\fill[qaone] (17,19) circle[radius=0.64pt];
\fill[qatwo] (17,20) circle[radius=0.64pt];
\fill[qatwo] (17,21) circle[radius=0.64pt];
\fill[qaone] (17,22) circle[radius=0.64pt];
\fill[qaone] (17,25) circle[radius=0.64pt];
\fill[qaone] (17,26) circle[radius=0.64pt];
\fill[qaone] (17,27) circle[radius=0.64pt];
\fill[qatwo] (18,3) circle[radius=0.64pt];
\fill[qatwo] (18,4) circle[radius=0.64pt];
\fill[qatwo] (18,5) circle[radius=0.64pt];
\fill[qaone] (18,6) circle[radius=0.64pt];
\fill[qaone] (18,7) circle[radius=0.64pt];
\fill[qaone] (18,8) circle[radius=0.64pt];
\fill[qatwo] (18,9) circle[radius=0.64pt];
\fill[qatwo] (18,10) circle[radius=0.64pt];
\fill[qatwo] (18,11) circle[radius=0.64pt];
\fill[qatwo] (18,12) circle[radius=0.64pt];
\fill[qatwo] (18,13) circle[radius=0.64pt];
\fill[qatwo] (18,14) circle[radius=0.64pt];
\fill[qatwo] (18,15) circle[radius=0.64pt];
\fill[qaone] (18,16) circle[radius=0.64pt];
\fill[qaone] (18,17) circle[radius=0.64pt];
\fill[qaone] (18,18) circle[radius=0.64pt];
\fill[qatwo] (18,19) circle[radius=0.64pt];
\fill[qatwo] (18,20) circle[radius=0.64pt];
\fill[qatwo] (18,21) circle[radius=0.64pt];
\fill[qathree] (18,24) circle[radius=0.64pt];
\fill[qathree] (18,25) circle[radius=0.64pt];
\fill[qaone] (18,26) circle[radius=0.64pt];
\fill[qaone] (18,27) circle[radius=0.64pt];
\fill[qatwo] (19,4) circle[radius=0.64pt];
\fill[qatwo] (19,5) circle[radius=0.64pt];
\fill[qatwo] (19,6) circle[radius=0.64pt];
\fill[qaone] (19,7) circle[radius=0.64pt];
\fill[qatwo] (19,8) circle[radius=0.64pt];
\fill[qatwo] (19,9) circle[radius=0.64pt];
\fill[qatwo] (19,10) circle[radius=0.64pt];
\fill[qatwo] (19,11) circle[radius=0.64pt];
\fill[qatwo] (19,12) circle[radius=0.64pt];
\fill[qatwo] (19,13) circle[radius=0.64pt];
\fill[qatwo] (19,14) circle[radius=0.64pt];
\fill[qatwo] (19,15) circle[radius=0.64pt];
\fill[qatwo] (19,16) circle[radius=0.64pt];
\fill[qaone] (19,17) circle[radius=0.64pt];
\fill[qatwo] (19,18) circle[radius=0.64pt];
\fill[qatwo] (19,19) circle[radius=0.64pt];
\fill[qatwo] (19,20) circle[radius=0.64pt];
\fill[qathree] (19,23) circle[radius=0.64pt];
\fill[qathree] (19,24) circle[radius=0.64pt];
\fill[qathree] (19,25) circle[radius=0.64pt];
\fill[qathree] (19,26) circle[radius=0.64pt];
\fill[qaone] (19,27) circle[radius=0.64pt];
\fill[qaone] (20,2) circle[radius=0.64pt];
\fill[qatwo] (20,5) circle[radius=0.64pt];
\fill[qatwo] (20,6) circle[radius=0.64pt];
\fill[qatwo] (20,7) circle[radius=0.64pt];
\fill[qatwo] (20,8) circle[radius=0.64pt];
\fill[qatwo] (20,9) circle[radius=0.64pt];
\fill[qatwo] (20,10) circle[radius=0.64pt];
\fill[qatwo] (20,11) circle[radius=0.64pt];
\fill[qaone] (20,12) circle[radius=0.64pt];
\fill[qatwo] (20,13) circle[radius=0.64pt];
\fill[qatwo] (20,14) circle[radius=0.64pt];
\fill[qatwo] (20,15) circle[radius=0.64pt];
\fill[qatwo] (20,16) circle[radius=0.64pt];
\fill[qatwo] (20,17) circle[radius=0.64pt];
\fill[qatwo] (20,18) circle[radius=0.64pt];
\fill[qatwo] (20,19) circle[radius=0.64pt];
\fill[qaone] (20,22) circle[radius=0.64pt];
\fill[qathree] (20,23) circle[radius=0.64pt];
\fill[qathree] (20,24) circle[radius=0.64pt];
\fill[qathree] (20,25) circle[radius=0.64pt];
\fill[qathree] (20,26) circle[radius=0.64pt];
\fill[qathree] (20,27) circle[radius=0.64pt];
\fill[qaone] (21,1) circle[radius=0.64pt];
\fill[qaone] (21,2) circle[radius=0.64pt];
\fill[qaone] (21,3) circle[radius=0.64pt];
\fill[qatwo] (21,6) circle[radius=0.64pt];
\fill[qatwo] (21,7) circle[radius=0.64pt];
\fill[qatwo] (21,8) circle[radius=0.64pt];
\fill[qatwo] (21,9) circle[radius=0.64pt];
\fill[qatwo] (21,10) circle[radius=0.64pt];
\fill[qaone] (21,11) circle[radius=0.64pt];
\fill[qaone] (21,12) circle[radius=0.64pt];
\fill[qaone] (21,13) circle[radius=0.64pt];
\fill[qatwo] (21,14) circle[radius=0.64pt];
\fill[qatwo] (21,15) circle[radius=0.64pt];
\fill[qatwo] (21,16) circle[radius=0.64pt];
\fill[qatwo] (21,17) circle[radius=0.64pt];
\fill[qatwo] (21,18) circle[radius=0.64pt];
\fill[qaone] (21,21) circle[radius=0.64pt];
\fill[qaone] (21,22) circle[radius=0.64pt];
\fill[qaone] (21,23) circle[radius=0.64pt];
\fill[qathree] (21,24) circle[radius=0.64pt];
\fill[qathree] (21,25) circle[radius=0.64pt];
\fill[qathree] (21,26) circle[radius=0.64pt];
\fill[qathree] (21,27) circle[radius=0.64pt];
\fill[qaone] (22,0) circle[radius=0.64pt];
\fill[qaone] (22,1) circle[radius=0.64pt];
\fill[qaone] (22,2) circle[radius=0.64pt];
\fill[qaone] (22,3) circle[radius=0.64pt];
\fill[qaone] (22,4) circle[radius=0.64pt];
\fill[qaone] (22,7) circle[radius=0.64pt];
\fill[qatwo] (22,8) circle[radius=0.64pt];
\fill[qatwo] (22,9) circle[radius=0.64pt];
\fill[qaone] (22,10) circle[radius=0.64pt];
\fill[qaone] (22,11) circle[radius=0.64pt];
\fill[qaone] (22,12) circle[radius=0.64pt];
\fill[qaone] (22,13) circle[radius=0.64pt];
\fill[qaone] (22,14) circle[radius=0.64pt];
\fill[qatwo] (22,15) circle[radius=0.64pt];
\fill[qatwo] (22,16) circle[radius=0.64pt];
\fill[qaone] (22,17) circle[radius=0.64pt];
\fill[qaone] (22,20) circle[radius=0.64pt];
\fill[qaone] (22,21) circle[radius=0.64pt];
\fill[qaone] (22,22) circle[radius=0.64pt];
\fill[qaone] (22,23) circle[radius=0.64pt];
\fill[qaone] (22,24) circle[radius=0.64pt];
\fill[qathree] (22,25) circle[radius=0.64pt];
\fill[qathree] (22,26) circle[radius=0.64pt];
\fill[qaone] (22,27) circle[radius=0.64pt];
\fill[qaone] (23,1) circle[radius=0.64pt];
\fill[qaone] (23,2) circle[radius=0.64pt];
\fill[qaone] (23,3) circle[radius=0.64pt];
\fill[qatwo] (23,8) circle[radius=0.64pt];
\fill[qatwo] (23,9) circle[radius=0.64pt];
\fill[qatwo] (23,10) circle[radius=0.64pt];
\fill[qaone] (23,11) circle[radius=0.64pt];
\fill[qaone] (23,12) circle[radius=0.64pt];
\fill[qaone] (23,13) circle[radius=0.64pt];
\fill[qatwo] (23,14) circle[radius=0.64pt];
\fill[qatwo] (23,15) circle[radius=0.64pt];
\fill[qatwo] (23,16) circle[radius=0.64pt];
\fill[qathree] (23,19) circle[radius=0.64pt];
\fill[qathree] (23,20) circle[radius=0.64pt];
\fill[qaone] (23,21) circle[radius=0.64pt];
\fill[qaone] (23,22) circle[radius=0.64pt];
\fill[qaone] (23,23) circle[radius=0.64pt];
\fill[qathree] (23,24) circle[radius=0.64pt];
\fill[qathree] (23,25) circle[radius=0.64pt];
\fill[qathree] (23,26) circle[radius=0.64pt];
\fill[qathree] (23,27) circle[radius=0.64pt];
\fill[qaone] (24,2) circle[radius=0.64pt];
\fill[qatwo] (24,9) circle[radius=0.64pt];
\fill[qatwo] (24,10) circle[radius=0.64pt];
\fill[qatwo] (24,11) circle[radius=0.64pt];
\fill[qaone] (24,12) circle[radius=0.64pt];
\fill[qatwo] (24,13) circle[radius=0.64pt];
\fill[qatwo] (24,14) circle[radius=0.64pt];
\fill[qatwo] (24,15) circle[radius=0.64pt];
\fill[qathree] (24,18) circle[radius=0.64pt];
\fill[qathree] (24,19) circle[radius=0.64pt];
\fill[qathree] (24,20) circle[radius=0.64pt];
\fill[qathree] (24,21) circle[radius=0.64pt];
\fill[qaone] (24,22) circle[radius=0.64pt];
\fill[qathree] (24,23) circle[radius=0.64pt];
\fill[qathree] (24,24) circle[radius=0.64pt];
\fill[qathree] (24,25) circle[radius=0.64pt];
\fill[qathree] (24,26) circle[radius=0.64pt];
\fill[qathree] (24,27) circle[radius=0.64pt];
\fill[qaone] (25,7) circle[radius=0.64pt];
\fill[qatwo] (25,10) circle[radius=0.64pt];
\fill[qatwo] (25,11) circle[radius=0.64pt];
\fill[qatwo] (25,12) circle[radius=0.64pt];
\fill[qatwo] (25,13) circle[radius=0.64pt];
\fill[qatwo] (25,14) circle[radius=0.64pt];
\fill[qaone] (25,17) circle[radius=0.64pt];
\fill[qathree] (25,18) circle[radius=0.64pt];
\fill[qathree] (25,19) circle[radius=0.64pt];
\fill[qathree] (25,20) circle[radius=0.64pt];
\fill[qathree] (25,21) circle[radius=0.64pt];
\fill[qathree] (25,22) circle[radius=0.64pt];
\fill[qathree] (25,23) circle[radius=0.64pt];
\fill[qathree] (25,24) circle[radius=0.64pt];
\fill[qathree] (25,25) circle[radius=0.64pt];
\fill[qathree] (25,26) circle[radius=0.64pt];
\fill[qaone] (25,27) circle[radius=0.64pt];
\fill[qaone] (26,6) circle[radius=0.64pt];
\fill[qaone] (26,7) circle[radius=0.64pt];
\fill[qaone] (26,8) circle[radius=0.64pt];
\fill[qatwo] (26,11) circle[radius=0.64pt];
\fill[qatwo] (26,12) circle[radius=0.64pt];
\fill[qatwo] (26,13) circle[radius=0.64pt];
\fill[qaone] (26,16) circle[radius=0.64pt];
\fill[qaone] (26,17) circle[radius=0.64pt];
\fill[qaone] (26,18) circle[radius=0.64pt];
\fill[qathree] (26,19) circle[radius=0.64pt];
\fill[qathree] (26,20) circle[radius=0.64pt];
\fill[qathree] (26,21) circle[radius=0.64pt];
\fill[qathree] (26,22) circle[radius=0.64pt];
\fill[qathree] (26,23) circle[radius=0.64pt];
\fill[qathree] (26,24) circle[radius=0.64pt];
\fill[qathree] (26,25) circle[radius=0.64pt];
\fill[qaone] (26,26) circle[radius=0.64pt];
\fill[qaone] (26,27) circle[radius=0.64pt];
\fill[qaone] (27,2) circle[radius=0.64pt];
\fill[qaone] (27,5) circle[radius=0.64pt];
\fill[qaone] (27,6) circle[radius=0.64pt];
\fill[qaone] (27,7) circle[radius=0.64pt];
\fill[qaone] (27,8) circle[radius=0.64pt];
\fill[qaone] (27,9) circle[radius=0.64pt];
\fill[qaone] (27,12) circle[radius=0.64pt];
\fill[qaone] (27,15) circle[radius=0.64pt];
\fill[qaone] (27,16) circle[radius=0.64pt];
\fill[qaone] (27,17) circle[radius=0.64pt];
\fill[qaone] (27,18) circle[radius=0.64pt];
\fill[qaone] (27,19) circle[radius=0.64pt];
\fill[qathree] (27,20) circle[radius=0.64pt];
\fill[qathree] (27,21) circle[radius=0.64pt];
\fill[qaone] (27,22) circle[radius=0.64pt];
\fill[qathree] (27,23) circle[radius=0.64pt];
\fill[qathree] (27,24) circle[radius=0.64pt];
\fill[qaone] (27,25) circle[radius=0.64pt];
\fill[qaone] (27,26) circle[radius=0.64pt];
\fill[qaone] (27,27) circle[radius=0.64pt];
\draw[->,semithick] (0,0) -- (28.5,0) node[right] {$m_0$};
\draw[->,semithick] (0,0) -- (0,28.5) node[above] {$m_1$};
\foreach \t in {0,9,18,27} {\draw (\t,0) -- (\t,-0.45) node[below] {\scriptsize\t};\draw (0,\t) -- (-0.45,\t) node[left] {\scriptsize\t};}
\end{scope}
\foreach \col/\lab/\pos in {qaone/1/0,qatwo/2/10,qathree/3/20} {\fill[\col!25] (\pos,-7) rectangle ++(3,2);\fill[\col] (\pos+1.5,-6) circle[radius=1.2pt];\node[anchor=west] at (\pos+4,-6) {$a=\lab$};}
\end{tikzpicture}
    \caption{Visualization of the values of $m=(m_0,m_1)$ for which either \eqref{eq:evendet1} or \eqref{eq:evendet2} holds for $n=6, p=5, k=1$. Shading shows the regions satisfying the inequalities, and colored dots are the exact nonnegative integer solutions. Pale gray dots satisfy neither condition.
    Where several exponents work, the smallest $a\geq1$ determines the color.
    All exponents are accounted for in this window: $a\geq4$ gives no solutions.}
    \label{figure}
\end{figure}

In the cases where the Hilbert series of the quasi-invariant polynomials in characteristic $p$ is different from in characteristic zero, we also provide an explicit Hilbert series in Corollary \ref{cor:degs}.

In Section \ref{sec:genfacts}, we state some general structural facts about quasi-invariant polynomials for dihedral groups. In Section \ref{sec:harmonic}, we discuss the space of $m$-harmonic polynomials, and give an explicit basis in terms of hypergeometric functions. In Section \ref{sec:reduction}, we introduce a reduction modulo $p$ map, and use the explicit description of $m$-harmonic polynomials to derive a necessary condition for the Hilbert series of the quasi-invariant polynomials in characteristic $p$ to be different from in characteristic zero. In Section \ref{sec:renxu}, we construct certain Frobenius-type quasi-invariant polynomials in characteristic $p$ and prove that they describe the lowest degree generators of $Q_m(D_n,\k_p)_{\rho_k}$ when the Hilbert series is different from in characteristic zero.

\section{Generalities on quasi-invariant polynomials for dihedral groups}\label{sec:genfacts}

In this section we prove some general structural facts about quasi-invariant polynomials for $W=D_n$ over any field $\k$ with $\char\k\nmid 2n$. These facts are all already known when $\char\k=0$ (see \cite{feigin}, \cite{eg}), but it will still be useful to consider this case in our context.

The quasi-invariance condition is homogeneous, so the quasi-invariant polynomials have an obvious grading by polynomial degree. The homogeneous degree $d$ polynomials in $\k[\h]$ are given by $r^dP(\theta)$, where $P(\theta)=\sum_{l=-d}^da_le^{li\theta}$ is a trigonometric polynomial such that $a_j=0$ if $d-j\equiv 1\pmod{2}$.

The generators of $\k[\h]^W$ are given by
\[e_2:=r^2,\quad e_n:=r^n(e^{ni\theta}+e^{-ni\theta}).\]

\subsection{One-dimensional representations}

Let
\[w:=\prod_{s\in\Sigma}\alpha_s(x),\quad w_0:=\prod_{j\equiv0\pmod{2}}\alpha_j(x),\quad w_1:=\prod_{j\equiv1\pmod{2}}\alpha_j(x),\quad w_m=\prod_{s\in\Sigma}\alpha_s(x)^{2m_s+1}.\]
We recall the description of the isotypic components $Q_m(D_n,\k)_\tau$ for $\dim\tau=1$.

\begin{proposition}[\cite{feigindihedral}, \cite{judy}]\label{prop:onedim}
    As a $\k[\h]^{D_n}$-module, we have
    \begin{enumerate}
        \item $Q_m(D_n,\k)_\triv$ is freely generated by $1$.
        \item $Q_m(D_n,\k)_\sgn$ is freely generated by $w_m$.
        \item $Q_m(D_n,\k)_{\rho_{\pm}}$ is freely generated by $w_1^{2m_1+1}$.
        \item $Q_m(D_n,\k)_{\rho_{\mp}}$ is freely generated by $w_0^{2m_0+1}$.
    \end{enumerate}
\end{proposition}

Thus the Hilbert series of the direct sum of $Q_m(D_n,\k)_\tau$ for $\dim\tau=1$ is
\[\frac{1+t^{2nm+n}}{(1-t^2)(1-t^n)}\]
if $n$ is odd and
\[\frac{1+t^{nm_0+n/2}+t^{nm_1+n/2}+{t^{n\mu+n}}}{(1-t^2)(1-t^n)}\]
if $n$ is even.

\subsection{Two-dimensional representations}

From the definition of $\rho_k$ it is clear that we have
\begin{lemma}\label{lem:dihedraldegs}
    Let
    \[f=r^d\sum_{l=-d}^da_le^{li\theta}\in\k[\h]\]
    be a homogeneous polynomial with $a_l=0$ when $d-l$ is odd. Then $f$ is in the $\rho_k$-isotypic component if and only if $a_l=0$ unless $l\equiv \pm k\pmod{n}$. 
\end{lemma}

Let $\rho_k$ be a two-dimensional representation of $D_n$ and let $Q_m(D_n,\k)_{\rho_k}^0$ be the $e^{2k\pi i/n}$-eigenspace of $t$ in $Q_m(D_n,\k)_{\rho_k}$. Then $Q_m(D_n,\k)_{\rho_k}^0$ is a $\k[\h]^{D_n}$-module, and we have an isomorphism of $\k[\h]^{D_n}$-modules
\[Q_m(D_n,\k)_{\rho_k}^0\otimes\rho_k\xto{\sim}Q_m(D_n,\k)_{\rho_k}\]
that sends $f\otimes v_k\mapsto f$. So to describe the $\k[\h]^{D_n}$-module structure of $Q_m(D_n,\k)_{\rho_k}$ it suffices to describe the $\k[\h]^{D_n}$-module structure of $Q_m(D_n,\k)_{\rho_k}^0$.

The rest of this section will be a generalization of some of the results in Section 3 of \cite{wang} to dihedral groups.

Note that for any two polynomials $f,g\in Q_m(D_n,\k)_{\rho_k}^0$, the polynomial $fs_0g-gs_0f$ lies in the isotypic component of $Q_m(D_n,\k)$ corresponding to $\wedge^2\rho_k\cong\sgn$. In particular, it is a $\k[\h]^{D_n}$-multiple of $w_m$ by Proposition \ref{prop:onedim}. Let us define a \emph{determinant} map
\[\det:Q_m(D_n,\k)_{\rho_k}^0\times Q_m(D_n,\k)_{\rho_k}^0\to\k[\h]^{D_n}\]
such that for $f,g\in Q_m(D_n,\k)_{\rho_k}^0$, we have
\[fs_0g-gs_0f=\det(f,g)w_m.\]

\begin{remark}
    If $\char\k=0$, then this is essentially the determinant considered in \cite[Lemma 5.18]{eg}.
\end{remark}

\begin{proposition}\label{prop:freeness}
\begin{enumerate}
    \item $Q_m(D_n,\k)_{\rho_k}^0$ is a free $\k[\h]^{D_n}$-module.\label{prop:freeness1}
    \item $Q_m(D_n,\k)_{\rho_k}^0$ has rank $2$ over $\k[\h]^{D_n}$. If $f,g$ are two free generators of $Q_m(D_n,\k)_{\rho_k}^0$, then we have $\det(f,g)\in\k^\times$, and thus $\deg f+\deg g=\deg w_m$.\label{prop:freeness2}
\end{enumerate}
\end{proposition}

\begin{proof}
    \ref{prop:freeness1}

    Note that by Quillen-Suslin, it suffices to show that $Q_m(D_n,\k)_{\rho_k}^0$ is projective. We will actually show that $Q_m(D_n,\k)$ is projective. Since $Q_m(D_n,\k)_{\rho_k}^0$ is a direct summand of $Q_m(D_n,\k)$, this immediately implies the result.
    
    Projectiveness is a local property, so it suffices to show that for all maximal ideals $\m\subset\k[\h]^{D_n}$, $Q_m(D_n,\k)_\m$ is a projective $\k[\h]^{D_n}_\m$-module. Since $\k[\h]^{D_n}_\m$ is a regular local ring of Krull dimension 2, by Auslander-Buchsbaum it suffices to show that the depth of $Q_m(D_n,\k)_\m$ is at least $2$. By \cite[Lemma 15.24.16]{stacks-project}, it suffices to show $Q_m(D_n,\k)$ is a reflexive module, and for this it is enough to show that $Q_m(D_n,\k)=\bigcap_\p Q_m(D_n,\k)_\p$ where the intersection is over all prime ideals $\p\subset\k[\h]^{D_n}$ of height $1$ by \cite[Lemma 15.24.18]{stacks-project}.

    Fix some $s\in\Sigma$, and consider the intersection $(\alpha_s(x))\cap\k[\h]^{D_n}$. Let $C$ be the conjugacy class of $s$. Any polynomial $e$ in the intersection must contain $\prod_{s'\in C}\alpha_{s'}$ as a factor. Then $e/\prod_{s'\in C}\alpha_{s'}$ is contained in a copy of the one-dimensional representation where $C$ acts by $-1$ and all other reflections act by $1$. But by Proposition \ref{prop:onedim}, the isotypic component of $\k[\h]=Q_0(D_n,\k)$ corresponding to this representation is generated by $\prod_{s'\in C}\alpha_{s'}$, so $e$ is divisible by $\prod_{s'\in C}\alpha_{s'}^2$. But the latter is invariant, so we have $(\alpha_s)\cap\k[\h]^{D_n}=(\prod_{s'\in C}\alpha_{s'}^2)$. In particular, there is a prime ideal of height 1 containing $(\alpha_s)\cap\k[\h]^{D_n}$, let it be $\p_s$.

    Then let $f\in\bigcap_\p Q_m(D_n,\k)_\p$. Then $f\in\bigcap_\p \k[\h]_\p=\k[\h]$. Since $f\in Q_m(D_n,\k)_{\p_s}$, there is some $e\not\in\p_s$ such that $f=f'/e$ for $f'\in Q_m(D_n,\k)$. Then we have $e(1-s)f=(1-s)f'$ is divisible by $\alpha_s(x)^{2m_s+1}$. But $e\not\in\p_s\supset(\alpha_s)\cap\k[\h]^{D_n}$, so this implies $(1-s)f$ is divisible by $\alpha_s(x)^{2m_s+1}$. Thus $f\in Q_m(D_n,\k)$.

    \ref{prop:freeness2}

    Consider the map
    \[\Phi:\k[\h]\to\bigoplus_{s\in\Sigma}\k[\h]/(\alpha_s(x)^{2m_s+1}),\quad f\mapsto(f\pmod{\alpha_s(x)^{2m_s+1}})_{s\in\Sigma}.\]
    Clearly, $\ker\Phi=(w_m)$. Note also that whether a polynomial is $m$-quasi-invariant only depends on its image under $\Phi$. More explicitly, for a given $s\in\Sigma$, let $\beta_s\in\h^*$ be a vector with $s\beta_s=\beta_s$. Then $\alpha_s,\beta_s$ form a basis of $\h^*$, and we have a basis of $\k[\h]/(\alpha_s(x)^{2m_s+1})$ given by $\alpha_s(x)^{d_1}\beta_s(x)^{d_2}$ with $d_1\leq 2m_s$. Let $Q_s$ be the subspace spanned by such basis vectors with $d_1$ even. Then for a polynomial $f$, $(1-s)f$ is divisible by $\alpha_s(x)^{2m_s+1}$ if and only if its image under $\Phi$ lies in $Q_s$. So we have
    \[Q_m(D_n,\k)=\Phi^{-1}\left(\bigoplus_{s\in\Sigma}Q_s\right).\]
    
    Next, we will show that $\Phi$ is surjective in degrees $d\geq\deg w_m$. Indeed, we have
    \[\dim((\k[\h]/(w_m))_d)=\dim(\k[\h]_d)-\dim(\k[\h]_{d-\deg w_m})=\deg w_m\]
    and
    \[\dim\left(\bigoplus_{s\in\Sigma}\k[\h]/(\alpha_s(x)^{2m_s+1})\right)_d=\sum_{s\in\Sigma}(\dim(\k[\h]_d)-\dim(\k[\h]_{d-(2m_s+1)}))=\sum_{s\in\Sigma}(2m_s+1)=\deg w_m.\]
    Thus since the map induced by $\Phi$ on $\k[\h]/(w_m)$ is injective, it must also be surjective in degrees at least $\deg w_m$ by dimension count.

    Then for any such $d$, the codimension of $Q_m(D_n,\k)_d$ in $\k[\h]_d$ is given by
    \[\dim(\k[\h]_d)-\dim(Q_m(D_n,\k)_d)=\dim((\k[\h]/(w_m))_d)-\dim((Q_m(D_n,\k)/(w_m))_d)=\deg(w_m)-\dim\left(\bigoplus_{s\in\Sigma}Q_s\right).\]
    In degree $d\geq\deg w_m\geq 2m_s+1$, we have $\dim(Q_s)=m_s+1$. So we have
    \[\deg(w_m)-\dim\left(\bigoplus_{s\in\Sigma}Q_s\right)=\sum_{s\in\Sigma}((2m_s+1)-(m_s+1))=\sum_{s\in\Sigma}m_s.\]
    This implies that we have
    \begin{equation}\label{eq:hilbbound}
        \lim_{t\to 1}(1-t)\left(\H(\k[\h])-\H(Q_m(D_n,\k))\right)=\sum_{s\in\Sigma}m_s.
    \end{equation}
    So the Gelfand-Kirillov dimension of the quotient (vector) space $\k[\h]/Q_m(D_n,\k)$ is $1$. In particular, this is only possible if $Q_m(D_n,\k)$ has as many free generators as $\k[\h]$ over $\k[\h]^{D_n}$. So each $Q_m(D_n,\k)_{\rho_k}^0$ has rank $2$.

    Let the free generators of $\k[\h]$ and $Q_m(D_n,\k)$ as $\k[\h]^{D_n}$-modules have degrees $d_j$ and $e_j$, respectively. Then we have
    \[\lim_{t\to 1}(1-t)\left(\H(\k[\h])-\H(Q_m(D_n,\k))\right)=\lim_{t\to 1}(1-t)\frac{\sum_j(t^{d_j}-t^{e_j})}{(1-t^2)(1-t^n)}=\frac{1}{2n}\sum_j(e_j-d_j).\]
    
    By Chevalley-Shephard-Todd, the coinvariant algebra $\k[\h]/\k[\h]^{D_n}_+\k[\h]$ is a Gorenstein ring, so we have $\sum_jd_j=|D_n|\frac{n}{2}=n^2$. For a two-dimensional representation $\rho_k$, let $f,g\in Q_m(D_n,\k)_{\rho_k}^0$ be free generators. Then if $\det(f,g)=0$, we obtain $\frac{f}{g}=\frac{s_0f}{s_0g}$ is fixed by $s_0$ and $t$, so it lies in $\k(\h)^{D_n}=\text{Frac}(\k[\h]^{D_n})$. But this implies a relation between $f$ and $g$, contradicting freeness. So we have $\det(f,g)\neq 0$. In particular, $\deg f+\deg g\geq\deg w_m$. 

    The degrees of the generators of $Q_m(D_n,\k)$ in isotypic components corresponding to one-dimensional representations are known via Proposition \ref{prop:onedim}. So the minimal possible value of $\sum_je_j$ is achieved when we have the equality case for each $\rho_k$, in which case we have $\sum_je_j=2n\frac{\deg w_m}{2}=2n\sum_{s\in\Sigma}m_s+n^2$. In this case, we have $\frac{1}{2n}\sum_j(e_j-d_j)=\sum_{s\in\Sigma}m_s$. But this is exactly what is required by \eqref{eq:hilbbound}. So we must have equality in the inequalities for each of the $\rho_k$. In particular, $\det(f,g)\neq 0$ is a degree zero element of $\k[\h]^{D_n}$, so must lie in $\k^\times$.
\end{proof}

In particular, this proposition implies

\begin{proposition}
    Let the minimal degree generator of $Q_m(D_n,\k)_{\rho_k}$ have degree $d$. Then we have
    \[\H(Q_m(D_n,\k)_{\rho_k})=\frac{2t^d+2t^{\deg w_m-d}}{(1-t^2)(1-t^n)}.\]
\end{proposition}

So to compute the Hilbert series of $Q_m(D_n,\k)$, it suffices to compute the lowest degree generators of $Q_m(D_n,\k)_{\rho_k}$. Moreover, we have

\begin{corollary}
    Let $\k[\h]^{D_n}_+$ denote the augmentation ideal of $\k[\h]^{D_n}$. Then $Q_m(D_n,\k)/\k[\h]^{D_n}_+$ is a Gorenstein ring.
\end{corollary}

\begin{proof}
    The top graded component of $Q_m(D_n,\k)/\k[\h]^{D_n}_+Q_m(D_n,\k)$ is spanned by $w_m$. Indeed, if there was another generator of $Q_m(D_n,\k)$ in degree $\deg w_m$, then it would have to lie in an isotypic component corresponding to a two-dimensional representation by Proposition \ref{prop:onedim}. But then Proposition \ref{prop:freeness} implies that the other generator of this isotypic component would lie in degree $0$, and the degree zero polynomials are constant.

    Thus the multiplication pairing on $Q_m(D_n,\k)/\k[\h]^{D_n}_+$ is the composition of multiplication and projection onto the $\sgn$ isotypic component. Let us denote this pairing by $\langle,\rangle$.

    It is easy to see from Proposition \ref{prop:onedim} that the multiplication pairing on $Q_m(D_n,\k)/\k[\h]^{D_n}_+$ is nondegenerate on each isotypic component corresponding to a one-dimensional representation. For a given two-dimensional representation $\rho_k$, let $f,g$ be the free generators of $Q_m(D_n,\k)_{\rho_k}^0$. Then up to scaling, $\{f,s_0f\}$ and $\{s_0g,g\}$ are dual bases of $D_nf,D_ng$ with respect to the multiplication pairing. Indeed, $\langle f,g\rangle$ and $\langle s_0f,s_0g\rangle$ are zero since the projection $\rho_k^{\otimes 2}\to\wedge^2\rho_k\cong\sgn$ sends the vectors $v_k\otimes v_k,s_0v_k\otimes s_0v_k$ to zero. On the other hand, applying $s_0$ yields $\langle f,s_0g\rangle=-\langle s_0f,g\rangle$. So we have
    \[\det(f,g)=\langle f,s_0g\rangle-\langle s_0f,g\rangle=2\langle f,s_0g\rangle=-2\langle s_0f,g\rangle\neq 0,\]
    and thus $\langle f,s_0g\rangle, \langle s_0f,g\rangle\neq 0$. So the multiplication pairing is nondegenerate.
\end{proof}

\begin{lemma}\label{lem:sympolys}
    A free generator $f\in Q_m(D_n,\k)_{\rho_k}^0$ is not divisible by $r^2$.
\end{lemma}

\begin{proof}
    If $f$ is divisible by $r^2$, then $f':=f/r^2$ is $m$-quasi-invariant since $(1-s)f'=(1-s)f/r^2$ is divisible by $\alpha_s(x)^{2m_s+1}$.
\end{proof}

\section{$m$-harmonic polynomials}\label{sec:harmonic}

Let us return to the setting where $W$ is any finite Coxeter group. For this section, let $\k$ be a field of characteristic 0 that contains $n$-th roots of unity. We will start by recalling results of \cite{feigin}, \cite{eg}, \cite{felder2003action}.

\begin{definition}
    The \textbf{Dunkl Laplacian} is the differential operator
    \begin{equation*}
        L=\Delta-\sum_{s\in\Sigma}\frac{2m_s}{\alpha_s(x)}\partial_{\alpha_s}
    \end{equation*}
    where $\Delta$ is the usual Laplacian.
\end{definition}

Note that this is the generalized Calogero-Moser operator after a gauge transformation (see \cite{feigin}). 

Let $N=\dim\h$, and let $e_1,\dots,e_N$ be the independent generators of $\k[\h]^W$ such that $e_1=\sum_{j=1}^Nx_j^2$ for an orthonormal basis $x_1,\dots,x_N$ of $\h$. Associated to each $e_j$ we have a quantum integral $L_j$ with leading term $e_j(\partial)$ with $L_1=L$. These $L_j$ are exactly the quantum integrals of the \emph{Calogero-Moser system}.

\begin{definition}
    A function $\psi$ is \textbf{$m$-harmonic} if
    \[L_j\psi=0\]
    for all $j$.
\end{definition}

Feigin and Veselov in \cite{feigin} reduced the Calogero-Moser problem to the study of $m$-harmonic functions. There, they proved that $m$-harmonic functions form a subspace of $Q_m(W,\k)$, which we will denote by $H_m$. The $m$-harmonic polynomials have the same degrees as the generators of $Q_m(W,\k)$ \cite{eg}, and their degrees were computed in \cite{felder2003action}. Moreover, in the case where $W$ is a dihedral group and $m$ is a constant multiplicity function, it was proved in \cite{feigin} that $Q_m(W,\k)$ is freely generated by $H_m$ as a $\k[\h]^W$-module. This was generalized to arbitrary multiplicity functions in \cite{feigindihedral}.

\subsection{$m$-harmonic polynomials for dihedral groups}

We will now compute the $m$-harmonic polynomials for the dihedral groups explicitly in terms of hypergeometric functions. Note that an alternative description of the $m$-harmonic polynomials without hypergeometric functions was computed in \cite{feigindihedral}, but it will be useful to the current paper to use the hypergeometric form.

Following the notation for $Q_m(D_n,\k)$, for an irreducible $D_n$-module $\rho$, we let $(H_m)_\rho$ be the isotypic component of $H_m$ corresponding to $\rho$, and if $\rho=\rho_k$ is a two-dimensional irreducible representation of $D_n$, then let $(H_m)_{\rho_k}^0$ be the $e^{2k\pi i/n}$-eigenspace of $t$ in $(H_m)_{\rho_k}$.

\begin{proposition}\label{prop:oddhyper}
    If $n$ is odd, a basis of $(H_m)_{\rho_k}^0$ is given by
    \begin{align*}
        h_{m,k,1}:=&e^{(mn+k)i\theta}{}_2F_1\left(-m,-m-\frac{k}{n};\frac{n-k}{n};e^{-2ni\theta}\right)r^{mn+k},\\
        h_{m,k,2}:=&e^{-(mn+n-k)i\theta}{}_2F_1\left(-m,-m-\frac{n-k}{n};\frac{k}{n};e^{2ni\theta}\right)r^{mn+n-k}.
    \end{align*}
\end{proposition}

\begin{proof}
    For ease of notation, assume that $\k=\C$. The quantum integrals $L_j$ and the polynomials $h_{m,k,1},h_{m,k,2}$ are defined over $\Q$, so this will not affect the validity of the proposition.
    
    By definition, all $m$-harmonic polynomials $f$ are solutions to $Lf=0$. By \cite{felder2003action}, we have
    \[\H(H_m)=1+\sum_{l=1}^{n-1}2t^{mn+l}+t^{n(2m+1)}.\]
    Note that an $m$-harmonic polynomial is not divisible by any positive degree symmetric polynomial. In particular, $r^2$ is such a polynomial, so an $m$-harmonic polynomial $f$ of degree $d$ has nonzero coefficient in either the top degree term $r^de^{di\theta}$ or the bottom degree term $r^de^{-di\theta}$. $f$ may be assumed to lie in some isotypic component, and by Lemma \ref{lem:dihedraldegs} this isotypic component must be $\rho_k$. So we see that $(H_m)_{\rho_k}$ consists of one copy of $\rho_k$ in degree $mn+k$ and one in degree $mn+n-k$.
    
    We will show that the given basis consists of all solutions to $Lf=0$ in these degrees with $f\in\k[\h]_{\rho_k}$. This will automatically imply that they are $m$-harmonic since there can be no other $m$-harmonic polynomials in these degrees.

    In \cite{feigin}, the following presentation for $L$ was computed in polar coordinates:
    \[L=\frac{1}{r^2}\partial_\theta^2-\frac{2mn}{r^2}\cot(n\theta)\partial_\theta+\partial_r^2+\frac{1-2mn}{r}\partial_r.\]
    Write $f=r^d\phi(\theta)$ where $d\in\{mn+k,mn+n-k\}$. Then the equation $Lf=0$ yields
    \begin{equation}\label{eq:oddhyper1}
        \partial_\theta^2\phi-2mn\cot(n\theta)\partial_\theta\phi+d(d-1)\phi+d(1-2mn)\phi=0.
    \end{equation}
    Recall that since $f$ is homogeneous of degree $d$, we may write $\phi$ as
    \[\phi(\theta)=\sum_{l=0}^da_le^{(d-2l)i\theta}.\]
    Since $f\in\k[\h]_{\rho_k}$, the only $a_l$ that are nonzero are the $l$ such that $d-2l\equiv k,n-k\pmod{n}$. If we choose $f$ to be in the $e^{2ki\pi/n}$-eigenspace of $t$, then $a_l=0$ for $d-2l\equiv n-k\pmod{n}$ as well, and if $f$ is in the $e^{-2ki\pi/n}$-eigenspace of $t$, then $a_l=0$ for $d-2l\equiv k\pmod{n}$. If $d=mn+k$, we choose $f$ to be in the $e^{2ki\pi/n}$-eigenspace of $t$, and if $d=mn+n-k$, we choose $f$ to be in the $e^{-2ki\pi/n}$-eigenspace of $t$. Then we obtain
    \[\phi(\theta)=\sum_{l=0}^{\frac{d}{n}}a_le^{(d-2nl)i\theta},\]
    so we may write $\phi=e^{di\theta}\psi(e^{-2ni\theta})$ for some polynomial $\psi$. Let $z:=e^{-2ni\theta}$. Then changing coordinates from $\theta$ to $z$, \eqref{eq:oddhyper1} becomes
    \begin{equation}\label{eq:oddhyper2}
        -4n^2z^2\partial_z^2\psi+\left(4n(d-n)z+4mn^2iz\cot(n\theta)\right)\partial_z\psi+\left(-d^2+d(d-2mn)-2mndi\cot(n\theta)\right)\psi=0.
    \end{equation}
    We have
    \[\cot(n\theta)=\frac{\frac{e^{ni\theta}+e^{-ni\theta}}{2}}{\frac{e^{ni\theta}-e^{-ni\theta}}{2i}}=i\frac{1+e^{-2ni\theta}}{1-e^{-2ni\theta}}=i\frac{1+z}{1-z},\]
    and after substituting into \eqref{eq:oddhyper2} and clearing denominators we obtain
    \begin{equation*}
        -4n^2z^2(1-z)\partial_z^2\psi+\left(4n(d-n)(1-z)-4mn^2(1+z)\right)z\partial_z\psi+\left(-2mnd(1-z)+2mnd(1+z)\right)\psi=0
    \end{equation*}
    which becomes
    \begin{equation*}
        z(1-z)\partial_z^2\psi+\left(m+1-\frac{d}{n}-z\left(1-m-\frac{d}{n}\right)\right)\partial_z\psi-\frac{md}{n}\psi=0
    \end{equation*}
    after simplifying. This is a hypergeometric differential equation, and has a unique polynomial solution given by
    \[\psi={}_2F_1\left(-m,-\frac{d}{n};m+1-\frac{d}{n};e^{-2ni\theta}\right)\]
    which is a Jacobi polynomial up to a change of variables. We obtain the $m$-harmonic polynomials by substituting in $d=mn+k,mn+n-k$.
\end{proof}

\begin{proposition}\label{prop:evenhyper}
    Let $n$ be even, and let 
    \[k_\mu=\begin{cases}
        k&\mu\text{ is even}\\
        \frac{n}{2}-k&\mu\text{ is odd}
    \end{cases}.\]
    Then a basis of $(H_m)_{\rho_k}$ is given by
    \begingroup
    \setlength{\thinmuskip}{2mu}
    \setlength{\medmuskip}{3mu}
    \small{\begin{align*}
        h_{m,k,1}:=& e^{-(n\mu/2+k_\mu)i\theta}(1-e^{ni\theta})^{\mu+2k_\mu/n}{}_2F_1\left(-\frac{\mu}{2}-\frac{k_\mu}{n},\frac{1+\delta}{2}-\frac{k_\mu}{n};1-\frac{2k_\mu}{n};-\frac{4e^{ni\theta}}{(1-e^{ni\theta})^2}\right)r^{n\mu/2+k_\mu},\\
        h_{m,k,2}:=& e^{(n\mu/2+n-k_\mu)i\theta}(1-e^{ni\theta})^{\mu+2(n-k_\mu)/n}{}_2F_1\left(\frac{k_\mu}{n}-1-\frac{\mu}{2},\frac{1+\delta}{2}-1+\frac{k_\mu}{n};\frac{2k_\mu}{n}-1;-\frac{4e^{ni\theta}}{(1e^{-ni\theta})^2}\right)r^{n\mu/2+n-k_\mu},
    \end{align*}}
    \endgroup
    $s_0h_{m,k,1}$, and $s_0h_{m,k,2}$.
\end{proposition}

Here, the expressions such as 
\[(1-e^{-ni\theta})^{\mu+2k_\mu/n}{}_2F_1\left(-\frac{\mu}{2}-\frac{k_\mu}{n},\frac{1+\delta}{2}-\frac{k_\mu}{n};1-\frac{2k_\mu}{n};-\frac{4e^{-ni\theta}}{(1-e^{-ni\theta})^2}\right)\]
should be viewed as formal power series in $e^{-ni\theta}$.

\begin{proof}
    Similarly to the odd case, assume $\k=\C$. From \cite{feigin} we obtain
    \[\frac{1}{\alpha_{s_j}(x)}\partial_{\alpha_j}=\frac{1}{r}\partial_r+\frac{1}{r^2}\cot\left(\theta-\frac{j\pi}{n}\right)\partial_\theta,\]
    so we have
    \begin{align*}
        L&=\Delta-2m_0\left(\frac{n}{2r}\partial_r+\frac{1}{r^2}\sum_{j\text{ even}}\cot\left(\theta-\frac{j\pi}{n}\right)\partial_\theta\right)-2m_1\left(\frac{n}{2r}\partial_r+\frac{1}{r^2}\sum_{j\text{ odd}}\cot\left(\theta-\frac{j\pi}{n}\right)\partial_\theta\right)\\
        &=\Delta-2m_0\left(\frac{n}{2r}\partial_r+\frac{n}{2r^2}\cot\left(\frac{n\theta}{2}\right)\partial_\theta\right)-2m_1\left(\frac{n}{2r}\partial_r+\frac{n}{2r^2}\cot\left(\frac{n\theta}{2}+\frac{\pi}{2}\right)\partial_\theta\right)\\
        &=\partial_r^2+\frac{1}{r^2}\partial_\theta^2+\frac{1}{r}\partial_r-\frac{n\mu}{r}\partial_r-\frac{n}{r^2}\left(m_0\cot\left(\frac{n\theta}{2}\right)-m_1\tan\left(\frac{n\theta}{2}\right)\right)\partial_\theta.
    \end{align*}
    By \cite{felder2003action}, we have
    \[\H(H_m)=1+t^{nm_0+n/2}+t^{nm_1+n/2}+\sum_{l=1}^{n/2-1}(2t^{n\mu/2+l}+2t^{n\mu/2+n-l})+t^{n(\mu+1)}.\]
    By Lemma \ref{lem:dihedraldegs}, the terms corresponding to $\rho_k$ are the terms in degree $\frac{n\mu}{2}+k_\mu,\frac{n\mu}{2}+n-k_\mu$.
    
    Let $d$ be one of $\frac{n\mu}{2}+k_\mu,\frac{n\mu}{2}+n-k_\mu$, and let $f\in H_m$ be a polynomial of degree $d$, and write $f=r^d\phi(\theta)$. Then the condition $Lf=0$ becomes
    \begin{equation}\label{eq:evenhyper}
        \partial_\theta^2\phi+n\left(m_1\tan\left(\frac{n\theta}{2}\right)-m_0\cot\left(\frac{n\theta}{2}\right)\right)\partial_\theta\phi+d(d-n\mu)\phi=0.
    \end{equation}
    Let $y=\sin^2\left(\frac{n\theta}{2}\right)$. Then writing $\phi$ in terms of $y$, \eqref{eq:evenhyper} becomes
    \[y(1-y)\partial_y^2\phi+\left(\frac{1}{2}-m_0-(1-\mu)y\right)\partial_y\phi+\frac{d}{n}\left(\frac{d}{n}-\mu\right)\phi=0.\]
    This is a hypergeometric differential equation, which has two linearly independent solutions near $y=\infty$. But note that the $\phi$ corresponding to the two linearly independent $m$-harmonic polynomials of degree $d$ are such solutions, so all solutions of the differential equation lead to $m$-harmonic polynomials. One such solution is
    \[y^{d/n}{}_2F_1\left(-\frac{d}{n},\frac{1}{2}+m_0-\frac{d}{n};1+\mu-\frac{2d}{n};y^{-1}\right),\]
    and the other is obtained by applying $s_0$ to this one. Noting that
    \[y^{d/n}=\left(\frac{e^{ni\theta/2}-e^{-ni\theta/2}}{2i}\right)^{2d/n}=\left(\frac{1}{2i}\right)^{2d/n}e^{di\theta}(1-e^{-ni\theta})^{2d/n}\]
    and
    \[y^{-1}=\csc^2\left(\frac{n\theta}{2}\right)=-\frac{4e^{ni\theta}}{(1-e^{ni\theta})^2},\]
    we obtain the proposition after rescaling by powers of $2i$ and substituting in values for $d$.
    
\end{proof}

\begin{remark}\label{rem:hyperid}
    Note that if we set $\delta=0$, then we are considering constant multiplicity functions $m$ for the even dihedral groups. In this case, replacing $n$ with an odd number in the proof of Proposition \ref{prop:evenhyper} yields a valid proof for the $m$-harmonic polynomials for odd dihedral groups. In particular, by comparison with Proposition \ref{prop:oddhyper}, we obtain identities such as
    \[{}_2F_1\left(-m,-m-\frac{k}{n};\frac{n-k}{n};e^{-2ni\theta}\right)=(1-e^{-ni\theta})^{2m+2k/n}{}_2F_1\left(-m-\frac{k}{n},\frac{1}{2}-\frac{k}{n};1-\frac{2k}{n};-\frac{4e^{-ni\theta}}{(1-e^{-ni\theta})^2}\right).\]
    These identities are examples of quadratic transformations of hypergeometric functions computed by Goursat.
\end{remark}

\section{Reduction modulo $p$}\label{sec:reduction}

For this section, fix a prime $p\nmid 2n$, and let $\k_0=\Q[\omega]$ be the $n$-th cyclotomic field. Let $R=\Z[\omega]$ be its ring of integers and let $\p\subset R$ be any prime ideal containing $p$. Then $\k_p:=R/\p$ is a field of characteristic $p$ containing $n$-th roots of unity, so we may use $\k_p$ as our base field for quasi-invariants in characteristic $p$. We will abuse notation and let $\h$ and $\rho_k$ refer to both the representations of $\k_0[D_n]$ and $\k_p[D_n]$. Then $R_\p[re^{i\theta},re^{-i\theta}]$ is an $R_\p$-lattice in $\k_0[\h]$, and the quotient map $R_\p[re^{i\theta},re^{-i\theta}]\to\k_p[\h]$ with kernel $\p R_\p[re^{i\theta},re^{-i\theta}]$ preserves quasi-invariants. We will refer to the restriction of this map to quasi-invariants as the \emph{reduction modulo $p$ map}
\[\pi_p:Q_m(D_n,\k_0)\cap R_\p[re^{i\theta},re^{-i\theta}]\to Q_m(D_n,\k_p).\]
Using this map, it follows from standard commutative algebra that $\H(Q_m(D_n,\k_p))\geq\H(Q_m(D_n,\k_0))$ in all degrees (see \cite[Proposition 2.2]{Ren_2020}).

For the rest of this section, let $k$ and $m$ be fixed. Note that all coefficients in the Taylor series expansions of $h_{m,k,1},h_{m,k,2}$ are rational, so there exist unique integers $\gamma_1,\gamma_2$ such that 
\[p^{\gamma_1}h_{m,k,1},p^{\gamma_2}h_{m,k,2}\in R[re^{i\theta},re^{-i\theta}]\setminus \p R[re^{i\theta},re^{-i\theta}].\]
$\gamma_1,\gamma_2$ are then the unique integers such that $\pi_p(p^{\gamma_1}h_{m,k,1}),\pi_p(p^{\gamma_2}h_{m,k,2})$ are both nonzero and well-defined. Moreover, $p^{\gamma_1}h_{m,k,1},p^{\gamma_2}h_{m,k,2}$ are the free generators of $Q_m(D_n,\k_0)_{\rho_k}^0$ by \cite{feigin}, \cite{feigindihedral}.

Using this, we obtain the following partial generalization of \cite[Theorem 4.1]{wang}:

\begin{proposition}\label{prop:reduction}
    If the Hilbert series of $Q_m(D_n,\k_p)_{\rho_k}$ differs from the Hilbert series of $Q_m(D_n,\k_0)_{\rho_k}$ then $\det(p^{\gamma_1}h_{m,k,1},p^{\gamma_2}h_{m,k,2})\in\p$.
\end{proposition}

\begin{proof}
    Since we have $\H(Q_m(D_n,\k_p)_{\rho_k}\geq\H(Q_m(D_n,\k_0)_{\rho_k}$ in all degrees, it is clear by Proposition \ref{prop:freeness} that the Hilbert series of the two differ if and only if the lowest degree generator of $Q_m(D_n,\k_p)_{\rho_k}^0$ has degree less than $\deg h_{m,k,1}$.

    If this is the case, then let $f$ be that generator. Then by Proposition \ref{prop:freeness}, $Q_m(D_n,\k_p)_{\rho_k}^0$ cannot have any other free generators of degree at most $\deg h_{m,k,2}$. In particular, $\pi_p(p^{\gamma_1}h_{m,k,1})$ and $\pi_p(p^{\gamma_2}h_{m,k,2})$ are both $\k_p[\h]^{D_n}$-multiples of $f$. Write $\pi_p(p^{\gamma_1}h_{m,k,1})=e_1f,\pi_p(p^{\gamma_2}h_{m,k,2})=e_2f$ so we have
    \begin{align*}
        \pi_p((p^{\gamma_1}h_{m,k,1})s_0(p^{\gamma_2}h_{m,k,2})-s_0(p^{\gamma_1}h_{m,k,1})(p^{\gamma_2}h_{m,k,2}))        &=e_1f(e_2s_0f)-e_1s_0f(e_2f)\\
        &=e_1e_2(fs_0f-fs_0f)\\
        &=0.
    \end{align*}
    So the reduction of $\det(p^{\gamma_1}h_{m,k,1},p^{\gamma_2}h_{m,k,2})$ modulo $p$ is zero and $\det(p^{\gamma_1}h_{m,k,1},p^{\gamma_2}h_{m,k,2})\in\p$.
\end{proof}

\begin{lemma}\label{lem:odddet}
    If $n$ is odd and $\det(p^{\gamma_1}h_{m,k,1},p^{\gamma_2}h_{m,k,2})\in\p$, then there exists some $a>0$ such that \eqref{eq:odddet} holds.
\end{lemma}
\begin{proof}
    The Taylor expansion of the hypergeometric function yields
    \begin{align*}
        h_{m,k,1}&=\left(\sum_{d=0}^m\frac{(-m)_d\left(-m-\frac{k}{n}\right)_d}{\left(\frac{n-k}{n}\right)_dd!}e^{(mn+k-2dn)i\theta}\right)r^{mn+k}\\
        s_0h_{m,k,2}&=\left(\sum_{d=0}^m\frac{(-m)_d\left(-m-\frac{n-k}{n}\right)_d}{\left(\frac{k}{n}\right)_dd!}e^{(mn+n-k-2dn)i\theta}\right)r^{mn+n-k}
    \end{align*}
    where $(q)_d=q(q+1)\dots(q+d-1)$ is the Pochhammer symbol. The highest degree term of $h_{m,k,1}s_0h_{m,k,2}$ in $e^{i\theta}$ is $(re^{i\theta})^{(2m+1)n}$ and the lowest degree term has degree $(2m+1)n-4mn=-(2m-1)n$ in $e^{i\theta}$. In particular, the highest degree term of $s_0(h_{m,k,1}s_0h_{m,k,2})$ has degree $(2m-1)n$, so the coefficient of $(re^{i\theta})^{(2m+1)n}$ in $h_{m,k,1}s_0h_{m,k,2}-s_0(h_{m,k,1}s_0h_{m,k,2})$ is $1$. The coefficient of $(re^{i\theta})^{(2m+1)n}$ in $w_m=(r^ne^{ni\theta}-r^ne^{-ni\theta})^{2m+1}$ is also $1$, so $\det(h_{m,k,1},h_{m,k,2})=1$. So $\det(p^{\gamma_1}h_{m,k,1},p^{\gamma_2}h_{m,k,2})\in\p$ if and only if $\gamma_1+\gamma_2>0$.

    We have
    \[\gamma_1=-\min_{0\leq d\leq m}v_p\left(\frac{(-m)_d\left(-m-\frac{k}{n}\right)_d}{\left(\frac{n-k}{n}\right)_dd!}\right),\quad \gamma_2=-\min_{0\leq d\leq m}v_p\left(\frac{(-m)_d\left(-m-\frac{n-k}{n}\right)_d}{\left(\frac{k}{n}\right)_dd!}\right).\]
    When $d=0$, both of these expressions have $p$-adic valuation $0$, so $\gamma_1,\gamma_2$ are nonnegative. So $\gamma_1+\gamma_2>0$ if and only if at least one of $\gamma_1,\gamma_2$ is positive.

    Let us first consider $\gamma_1$. Observe that if for some $d\leq m$ we have $v_p\left(\frac{(-m)_d\left(-m-\frac{k}{n}\right)_d}{\left(\frac{n-k}{n}\right)_dd!}\right)<0$, then after expanding out the Pochhammer symbols and factorial, there must be some $a>0$ such that there are more terms in the denominator that are divisible by $p^a$ than terms in the numerator divisible by $p^a$. Observe that if $q$ is a rational number whose denominator is not divisible by $p$, then the number of terms in the expansion of $(q)_d$ that are divisible by $p^a$ is exactly the number of nonnegative integers less than $d$ that are congruent to $-q$ mod $p^a$. This number always increases by $1$ when $d$ increases by $p^a$, and since $d!=(1)_d$ can also be expressed as a Pochhammer symbol, we obtain that the difference between the number of terms in $(-m)_d\left(-m-\frac{k}{n}\right)_d$ and $\left(\frac{n-k}{n}\right)_dd!$ that are divisible by $p^a$ is periodic in $d$ with period $p^a$. So we may assume $d<p^a$. Then $d!$ clearly has no terms divisible by $p^a$, and $(-m)_d$, $\left(-m-\frac{k}{n}\right)_d$, and $\left(\frac{n-k}{n}\right)_d$ each have at most one term divisible by $p^a$. So we require there to be some $d$ such that such a term exists in $\left(\frac{n-k}{n}\right)_d$ but not in the other two Pochhammer symbols. This happens exactly when $-\frac{n-k}{n}\mod{p^a}=\frac{k}{n}-1\mod{p^a}$ is less than both $m\mod{p^a}$ and $m+\frac{k}{n}\mod{p^a}$. The latter occurs exactly when $(m\mod{p^a})+(\frac{k}{n}\mod{p^a})<p^a$, which happens when $m\mod{p^a}<-\frac{k}{n}\mod{p^a}$.

    So if $\gamma_1$ is positive, then we have
    \[\frac{k}{n}-1\mod{p^a}<m\mod{p^a}<-\frac{k}{n}\mod{p^a}.\]
    Similarly, if $\gamma_2$ is positive, then
    \[-\frac{k}{n}\mod{p^a}<m\mod{p^a}<\frac{k}{n}-1\mod{p^a}.\]
    If $p^a|k$, then $\frac{k}{n}\mod{p^a}=0$, so $\gamma_1$ is never positive, and $\gamma_2$ is positive when $0<m\mod{p^a}<p^a-1$. We have $\Fr^a(k)=\frac{k}{p^a}$, so for integer $m$, this is exactly equivalent to \eqref{eq:odddet}.
    
    Now, assume $p^a$ does not divide $k$. If $\Fr^a(k)=kp^{-a}\mod{n}$, then we have $\frac{k}{n}\mod{p^a}=\frac{(n-\Fr^a(k))p^a+k}{n}$. In particular, we have $\frac{k}{n}\mod{p^a}>-\frac{k}{n}\mod{p^a}$, so $\gamma_1$ is never positive, and if $\gamma_2$ is positive, then
    \[\frac{\Fr^a(k)p^a-k}{n}<m\mod{p^a}<\frac{(n-\Fr^a(k))p^a-(n-k)}{n},\]
    which for integer $m$ is equivalent to \eqref{eq:odddet}. If $\Fr^a(k)=-kp^{-a}$, then similarly to the previous case, $\gamma_2$ is never positive and if $\gamma_1$ is positive, then $m$ satisfies \eqref{eq:odddet}.
\end{proof}

\begin{lemma}\label{lem:evendet}
    If $n$ is even and $\det(p^{\gamma_1}h_{m,k,1},p^{\gamma_2}h_{m,k,2})\in\p$, then there exists some $a>0$ such that \eqref{eq:evendet1} or \eqref{eq:evendet2} holds.
\end{lemma}

\begin{proof}
    The $t$-action on $h_{m,k,1}$ combined with degree constraints imply that the highest degree term of $h_{m,k,1}$ in $e^{i\theta}$ has degree $\frac{n\mu}{2}+k_\mu$ and the lowest degree term has degree at least $-\frac{n\mu}{2}+k_\mu$. Moreover, Proposition \ref{prop:evenhyper} implies that the highest degree term has coefficient $1$. Similarly, the lowest degree term of $h_{m,k,2}$ has degree $-\frac{n\mu}{2}-n+k\mu$ with coefficient $1$ and the highest degree term has degree at most $-\frac{n\mu}{2}+n-k_\mu$. By the same argument as in the previous lemma, we obtain $\det(p^{\gamma_1}h_{m,k,1},p^{\gamma_2}h_{m,k,2})\in\p$ if and only if $\gamma_1+\gamma_2>0$. Since $h_{m,k,1},h_{m,k,2}$ both have a coefficient equal to $1$, we have $\gamma_1,\gamma_2\geq 0$. 

    If $\gamma_1>0$, then by Proposition \ref{prop:evenhyper} there must be some coefficient in the hypergeometric function
    \[{}_2F_1\left(-\frac{\mu}{2}-\frac{k_\mu}{n},\frac{1+\delta}{2}-\frac{k_\mu}{n};1-\frac{2k_\mu}{n};-\frac{4e^{-ni\theta}}{(1-e^{-ni\theta})^2}\right)\]
    whose $p$-adic valuation is less than $0$, which can only occur if there is some $d$ such that
    \[v_p\left(\frac{\left(-\frac{\mu}{2}-\frac{k_\mu}{n}\right)_d\left(\frac{1+\delta}{2}-\frac{k_\mu}{n}\right)_d}{\left(1-\frac{2k_\mu}{n}\right)_dd!}\right)<0. \]
    By the same argument as in Lemma \ref{lem:odddet}, this occurs exactly when there is some $a$ such that $\frac{2k_\mu}{n}-1\mod{p^a}$ is less than both $\frac{\mu}{2}+\frac{k_\mu}{n}\mod{p^a}$ and $\frac{k_\mu}{n}-\frac{1+\delta}{2}\mod{p^a}$.

    If $p^a|k_\mu$, then we have $\frac{2k_\mu}{n}-1\mod{p^a}=p^a-1$, so we never have $\gamma_1>0$. Otherwise, note that depending on whether $p^a\Fr^a(k_\mu)$ is congruent to $k_\mu$ or $-k_\mu$ modulo $n$, we have that $\frac{k_\mu}{n}\mod p^a=\frac{(n-\Fr^a(k_\mu))p^a+k_\mu}{n}$ or $\frac{\Fr^a(k_\mu)p^a+k_\mu}{n}$, respectively. The former is always greater than $\frac{p^a}{2}$, and if the latter is greater than $\frac{p^a}{2}$ then we have $k_\mu>(\frac{n}{2}-\Fr^a(k_\mu))p^a$. But note that we have $(\frac{n}{2}-\Fr^a(k_\mu))p^a>0$, and $(\frac{n}{2}-\Fr^a(k_\mu))p^a\equiv\frac{n}{2}+k_\mu\pmod{n}$, so $(\frac{n}{2}-\Fr^a(k_\mu))p^a>\frac{n}{2}>k_\mu$ and this cannot occur. So we conclude that $\frac{\Fr^a(k_\mu)p^a+k_\mu}{n}<\frac{p^a}{2}$.

    Let us now do casework on whether $\frac{k_\mu}{n}\mod p^a$ is greater or less than $\frac{p^a}{2}$ and whether $\mu$ and $\delta$ are even or odd. If $\frac{k_\mu}{n}\mod p^a<\frac{p^a}{2}$, then $\frac{2k_\mu}{n}\mod p^a>\frac{k_\mu}{n}\mod p^a$, so $\frac{2k_\mu}{n}-1\mod p^a<\frac{\mu}{2}+\frac{k_\mu}{n}\mod p^a$ if $\frac{\mu}{2}\mod{p^a}>\frac{k_\mu}{n}-1\mod{p^a}$ and $(\frac{\mu}{2}\mod{p^a})+(\frac{k_\mu}{n}\mod{p^a})<p^a$. Writing $\frac{k_\mu}{n}\mod{p^a}=\frac{\Fr^a(k_\mu)p^a+k_\mu}{n}$, this happens when
    \begin{equation}\label{eq:mu1}
        \frac{\Fr^a(k_\mu)p^a-(n-k_\mu)}{n}<\frac{\mu}{2}\mod{p^a}<\frac{(n-\Fr^a(k_\mu))p^a-k_\mu}{n}.
    \end{equation}
    Similarly, we have $\frac{2k_\mu}{n}-1\mod{p^a}<\frac{k_\mu}{n}-\frac{1+\delta}{2}\mod{p^a}$ when
    \begin{equation}\label{eq:delta1}
        \frac{\Fr^a(k_\mu)p^a-(n-k_\mu)}{n}<-\frac{1+\delta}{2}\mod{p^a}<\frac{(n-\Fr^a(k_\mu))p^a-k_\mu}{n}.
    \end{equation}
    When $\mu,\delta$ are even, then we obtain $\frac{\mu}{2}\mod{p^a}=p^a\left\{\frac{\mu}{2p^a}\right\}$ and $-\frac{1+\delta}{2}\mod{p^a}=p^a\left\{-\frac{1+\delta+p^a}{2p^a}\right\}$, so \eqref{eq:mu1} and \eqref{eq:delta1} are equivalent to \eqref{eq:evendet1} for even integers $\mu,\delta$. Similarly, if $\mu,\delta$ are odd then we have $\frac{\mu}{2}\mod{p^a}=p^a\left\{\frac{\mu+p^a}{2p^a}\right\}$ and $-\frac{1+\delta}{2}\mod{p^a}=p^a\left\{-\frac{1+\delta}{2p^a}\right\}$. Note that in this case we have $k_\mu=\frac{n}{2}-k$, and it is clear from \eqref{eq:frob} that $\Fr^a(\frac{n}{2}-k)=\frac{n}{2}-\Fr^a(k)$, so \eqref{eq:mu1} and \eqref{eq:delta1} are equivalent to \eqref{eq:evendet2} for odd integers $\mu,\delta$.

    If $\frac{k_\mu}{n}\mod{p^a}>\frac{p^a}{2}$, then $\frac{2k_\mu}{n}\mod{p^a}<\frac{k_\mu}{n}\mod{p^a}$, so $\frac{2k_\mu}{n}-1\mod p^a<\frac{\mu}{2}+\frac{k_\mu}{n}\mod p^a$ if either $(\frac{\mu}{2}\mod{p^a})+(\frac{k_\mu}{n}\mod{p^a})<p^a$ or $\frac{\mu}{2}\mod{p^a}>\frac{k_\mu}{n}-1\mod{p^a}$. This happens when
    \begin{equation*}
        \frac{(n-\Fr^a(k_\mu))p^a-(n-k_\mu)}{n}<\frac{\mu}{2}\mod{p^a}\quad\text{or}\quad\frac{\mu}{2}\mod{p^a}<\frac{\Fr^a(k_\mu)p^a-k_\mu}{n}.
    \end{equation*}
    Similarly, we have $\frac{2k_\mu}{n}-1\mod{p^a}<\frac{k_\mu}{n}-\frac{1+\delta}{2}\mod{p^a}$ when
    \begin{equation*}
        \frac{(n-\Fr^a(k_\mu))p^a-(n-k_\mu)}{n}<-\frac{1+\delta}{2}\mod{p^a}\quad\text{or}\quad-\frac{1+\delta}{2}\mod{p^a}<\frac{\Fr^a(k_\mu)p^a-k_\mu}{n}.
    \end{equation*}
    When $\mu,\delta$ are even, these are equivalent to \eqref{eq:evendet2} and when $\mu,\delta$ are odd, they are equivalent to \eqref{eq:evendet1}.
    
    If $\gamma_2>0$, then there is some $d$ such that
    \[v_p\left(\frac{\left(\frac{k_\mu}{n}-1-\frac{\mu}{2}\right)_d\left(\frac{1+\delta}{2}-1+\frac{k_\mu}{n}\right)_d}{\left(\frac{2k_\mu}{n}-1\right)_dd!}\right)<0.\]
    If $p^a$ does not divide $k_\mu$, the same argument as for $\gamma_1$ implies that $\gamma_2>0$ only when \eqref{eq:evendet1} or \eqref{eq:evendet2} holds. When $p^a$ divides $k_\mu$, we can have $\gamma_2>0$ only when there is some $a$ such that $\frac{\mu}{2}+1\mod{p^a}>1$ and $\frac{1-\delta}{2}\mod{p^a}>1$. We have $\Fr^a(k_\mu)=\frac{k_\mu}{p^a}<\frac{n}{2p^a}$, so when $\mu$ is even, this is equivalent to \eqref{eq:evendet1} for even integers $\mu,\delta$ and when $\mu$ is odd, this is equivalent to \eqref{eq:evendet2} for odd integers $\mu,\delta$.
\end{proof}

When $\delta=0$, so that $m=\frac{\mu}{2}$ is a constant function, \eqref{eq:evendet1} reduces to \eqref{eq:odddet} and \eqref{eq:evendet2} is never satisfied. Therefore, using the ideas from Remark \ref{rem:hyperid}, this case yields an alternate proof to Lemma \ref{lem:odddet}.

\section{Frobenius-type quasi-invariants}\label{sec:renxu}

In this section we recall a construction of Ren-Xu in \cite{Ren_2020} to produce a class of quasi-invariant polynomials in characteristic $p$. For dihedral groups, when $\H(Q_m(D_n,\k_p))\neq\H(Q_m(D_n,\k_0))$, the generators of $Q_m(D_n,\k_p)_{\rho_k}$ for $1\leq k<\frac{n}{2}$ are obtained via this construction.

\begin{lemma}\label{lem:frob}
    Let $\k$ be a field of characteristic $p$ and let $m',b:\Sigma\to\Z_{\geq 0}$ be functions invariant on conjugacy classes and let $a\in\Z_{>0}$. Let $m=\frac{p^a(2m'+1)-1}{2}+b$. Then for any $f\in Q_{m'}(W,\k)$, $f^{p^a}\prod_{s\in\Sigma}\alpha_s(x)^{2b_s}\in Q_{m}(W,\k)$.
\end{lemma}

See \cite[Theorem 2.6]{Ren_2020} for the proof of this proposition when $W$ is a symmetric group; the same proof holds for general $W$ as well.

For even $n$, we also need the following modified version of the above construction.

\begin{lemma}\label{lem:twistedfrob}
    Let $n$ be even, let $\k$ be a field of characteristic $p$ and let $m',b:\Sigma\to\Z_{\geq 0}$ be functions invariant on conjugacy classes and let $a\in\Z_{>0}$. Let $j\in\{0,1\}$, and let $m_{1-j}=\frac{p^a(2m'_{1-j}+1)-1}{2}+b_{1-j}$, $m_j=b_j$. Then for any $f\in Q_{m'}(W,\k)$, $f^{p^a}w_j\prod_{s\in\Sigma}\alpha_s(x)^{2b_s}\in Q_{m}(W,\k)$.
\end{lemma}
\begin{proof}
    Let $g=f^{p^a}w_j\prod_{s\in\Sigma}\alpha_s(x)^{2b_s}=f^{p^a}w_jw_0^{2b_0}w_1^{2b_1}$.
    
    First, let $s\in\Sigma$ be such that $\alpha_s(x)$ is a factor of $w_j$. Then note that $g,sg$ are both divisible by $w_j^{2b_j+1}$. So $(1-s)g$ is divisible by $\alpha_s(x)^{2m_j+1}$. Next, let $s$ be such that $\alpha_s(x)$ is not a factor of $w_j$. Then $w_jw_0^{2b_0}w_1^{2b_1}$ is $s$-invariant. In particular, we have $(1-s)g=w_jw_0^{2b_0}w_1^{2b_1}(1-s)f^{p^a}$ is divisible by $\alpha_s(x)^{2m_{1-j}+1}$ by Lemma \ref{lem:frob}.
\end{proof}
We will refer to quasi-invariants obtained via the construction in these lemmas as \emph{Frobenius-type quasi-invariants}.

\begin{proposition}\label{prop:renxu}
    If $n$ is odd and $m,k$ are such that there exists some $a$ such that \eqref{eq:odddet} holds, then the Hilbert series of $Q_m(D_n,\k_p)_{\rho_k}$ is different from the Hilbert series of $Q_m(D_n,\k_0)_{\rho_k}$.
    
    If $n$ is even and $\mu,\delta,k$ are such that there exists some $a$ such that \eqref{eq:evendet1} or \eqref{eq:evendet2} holds, then the Hilbert series of $Q_m(D_n,\k_p)_{\rho_k}$ is different from the Hilbert series of $Q_m(D_n,\k_0)_{\rho_k}$.
\end{proposition}

\begin{proof}
    We will just prove the even case; the proof of the odd case is obtained by setting $\delta=0,m=\frac{\mu}{2}$ in the proof of the even case.

    Recall that the minimal degree generator of $Q_m(D_n,\k_0)_{\rho_k}$ has degree $\frac{n\mu}{2}+k_\mu$, so it suffices to show that when \eqref{eq:evendet1} or \eqref{eq:evendet2} holds, there is a Frobenius-type quasi-invariant in $Q_m(D_n,\k_p)_{\rho_k}$ of degree less than $\frac{n\mu}{2}+k_\mu$.

    First, assume that there is some $a$ such that \eqref{eq:evendet1} holds. Then let
    \[m_0'=\left\lfloor\frac{\mu}{2p^a}\right\rfloor+\left\lfloor\frac{\delta}{2p^a}+\frac{1}{2}\right\rfloor,\quad m_1'=\left\lfloor\frac{\mu}{2p^a}\right\rfloor-\left\lfloor\frac{\delta}{2p^a}+\frac{1}{2}\right\rfloor,\quad b=\max\left(0,m-\frac{p^a(2m'+1)-1}{2}\right).\]
    Let us first assume that $m_0',m_1'\geq 0$. Note that $m'_0+m'_1$ is even, so by the reduction modulo $p$ map, there is some $f\in Q_{m'}(D_n,\k_p)_{\rho_{\Fr^a(k)}}$ of degree $\frac{n(m'_0+m'_1)}{2}+\Fr^a(k)$. Then by Lemma \ref{lem:frob}, $f^{p^a}w_0^{2b_0}w_1^{2b_1}$ is an element of $Q_m(D_n,\k_p)$, and it lies in the isotypic component corresponding to the representation $(\Fr^a)^*(\rho_{\Fr^a(k)})\cong\rho_k$. Its degree is equal to
    \begin{equation}\label{eq:renxuexp}
        p^a\left(\frac{n(m'_0+m'_1)}{2}+\Fr^a(k)\right)+n(b_0+b_1)=np^a\left\lfloor\frac{\mu}{2p^a}\right\rfloor+p^a\Fr^a(k)+n(b_0+b_1).
    \end{equation}
    We wish for this to be less than $\frac{n\mu}{2}+k_\mu$. Note that we have
    \[\frac{n\mu}{2}=np^a\frac{\mu}{2p^a}=np^a\left\lfloor\frac{\mu}{2p^a}\right\rfloor+np^a\left\{\frac{\mu}{2p^a}\right\},\]
    so the inequality we wish to show is
    \begin{equation}\label{eq:renxudegs}
        np^a\left\{\frac{\mu}{2p^a}\right\}+k_\mu>p^a\Fr^a(k)+n(b_0+b_1).
    \end{equation}
    If $b_0=b_1=0$, then this follows immediately from \eqref{eq:evendet1}.
    If $b_0,b_1>0$, then we have
    \[n(b_0+b_1)=n\left(\mu-p^a(m_0'+m_1')-(p^a-1)\right)=n\mu-2np^a\left\lfloor\frac{\mu}{2p^a}\right\rfloor-n(p^a-1)=2np^a\left\{\frac{\mu}{2p^a}\right\}-n(p^a-1)\]
    and \eqref{eq:renxudegs} becomes
    \[(n-\Fr^a(k))p^a-(n-k_\mu)>np^a\left\{\frac{\mu}{2p^a}\right\},\]
    which follows from \eqref{eq:evendet1}.
    If $b_0>0,b_1=0$, then 
    \[nb_0=nm_0-np^am_0'-\frac{n(p^a-1)}{2}=nm_0-np^a\left\lfloor\frac{\mu}{2p^a}\right\rfloor-np^a\left\lfloor\frac{\delta}{2p^a}+\frac{1}{2}\right\rfloor-\frac{n(p^a-1)}{2}.\]
    Note that
    \[\left\lfloor\frac{\mu}{2p^a}\right\rfloor+\left\lfloor\frac{\delta}{2p^a}+\frac{1}{2}\right\rfloor=\frac{\mu}{2p^a}+\frac{\delta}{2p^a}+\frac{1}{2}-\left\{\frac{\mu}{2p^a}\right\}-\left\{\frac{\delta}{2p^a}+\frac{1}{2}\right\}=\frac{m_0}{p^a}+\frac{1}{2}-\left\{\frac{\mu}{2p^a}\right\}-\left\{\frac{\delta}{2p^a}+\frac{1}{2}\right\},\]
    so we have
    \[nb_0=nm_0-nm_0-\frac{np^a}{2}+np^a\left\{\frac{\mu}{2p^a}\right\}+np^a\left\{\frac{\delta}{2p^a}+\frac{1}{2}\right\}-\frac{n(p^a-1)}{2}=np^a\left\{\frac{\mu}{2p^a}\right\}+np^a\left\{\frac{\delta}{2p^a}+\frac{1}{2}\right\}-np^a+\frac{n}{2}\]
    and \eqref{eq:renxudegs} becomes
    \[(n-\Fr^a(k))p^a-\left(\frac{n}{2}-k_\mu\right)>np^a\left\{\frac{\delta}{2p^a}+\frac{1}{2}\right\}\]
    which follows from \eqref{eq:evendet1}. The case $b_0=0,b_1>0$ is similar.

    Now, consider when $m_0'$ or $m_1'$ is negative. Then note that $|\delta|\leq\mu$ since $m_0,m_1\geq 0$, so we have
    \[\left\lfloor\frac{\mu}{2p^a}\right\rfloor+1\geq\left|\left\lfloor\frac{\delta}{2p^a}+\frac{1}{2}\right\rfloor\right|,\]
    with the $+1$ correction coming from the floor functions. In particular, only one of $m_0'$ or $m_1'$ is negative, and if it is negative, it must be equal to $-1$. Without loss of generality, let us assume that $m_0'=-1$. Then set $m''_0=0,m''_1=m'_1$, and let $b'_0=m_0,b'_1=b_1$. Then $m_0''+m_1''$ is odd, so there is some $f\in Q_{m''}(D_n,\k_p)_{\rho_{n/2-\Fr^a(k)}}$ of degree $\frac{nm_1''}{2}+(\frac{n}{2}-(\frac{n}{2}-\Fr^a(k)))=\frac{nm_1''}{2}+\Fr^a(k)$. By Lemma \ref{lem:twistedfrob}, $f^{p^a}w_0^{2b'_0+1}w_1^{2b'_1}$ is an element of $Q_m(D_n,\k_p)$, and it lies in the isotypic component corresponding to $(\Fr^a)^*(\rho_{n/2-\Fr^a(k)})\otimes\rho_{\mp}\cong\rho_{n/2-k}\otimes\rho_\mp\cong\rho_k$. Its degree is equal to
    \[p^a\left(\frac{nm_1''}{2}+\Fr^a(k)\right)+\frac{n}{2}+n(b'_0+b'_1)=p^a\left(\frac{n(m_0'+m_1')}{2}+\Fr^a(k)\right)+\frac{n(p^a+1)}{2}+nm_0+nb_1.\]
    Since $m_0'=-1$, we have $b_0=m_0+\frac{p^a+1}{2}$, so this is exactly the same expression as \eqref{eq:renxuexp}. So the same analysis implies that the degree of $f^{p^a}w_0^{2b'_0+1}w_1^{2b'_1}$ is less than $\frac{n\mu}{2}+k_\mu$.

    If there is some $a$ such that \eqref{eq:evendet2} holds, then let
    \[m_0'=\left\lfloor\frac{\mu}{2p^a}+\frac{1}{2}\right\rfloor+\left\lfloor\frac{\delta}{2p^a}\right\rfloor,\quad m_1'=\left\lfloor\frac{\mu}{2p^a}+\frac{1}{2}\right\rfloor-\left\lfloor\frac{\delta}{2p^a}\right\rfloor-1,\quad b=\max\left(0,m-\frac{p^a(2m'+1)-1}{2}\right).\]
    
    If $m_0',m_1'\geq 0$, then $m_0'+m_1'$ is odd, so there is some $f\in Q_{m'}(D_n,\k_p)_{\rho_{\Fr^a(k)}}$ of degree $\frac{n(m_0'+m_1')}{2}+(\frac{n}{2}-\Fr^a(k))$. Then $f^{p^a}w_0^{2b_0}w_1^{2b_1}\in Q_m(D_n,\k_p)_{\rho_k}$, and a similar argument as in the previous case implies that it has degree less than $\frac{n\mu}{2}+k_\mu$.

    If $m_0'$ or $m_1'$ is negative, then we again have that only one of them can be negative, and then it must be $-1$. Assume it is $m_0'$, then define $m'',b'$ the same way as in the previous case. The same analysis shows that for $f\in Q_{m''}(D_n,\k_p)_{\rho_{n/2-\Fr^a(k)}}$ of degree $\frac{nm_1''}{2}+(\frac{n}{2}-\Fr^a(k))$, $f^{p^a}w_0^{2b'_0+1}w_1^{2b'_1}$ is an element of $Q_m(D_n,\k_p)_{\rho_k}$ of degree less than $\frac{n\mu}{2}+k_\mu$.
\end{proof}

Combining Proposition \ref{prop:reduction}, Lemmas \ref{lem:odddet} and \ref{lem:evendet}, and Proposition \ref{prop:renxu} yields Theorem \ref{thm:condition}. Moreover, we can now compute the lowest degree generator of $Q_m(D_n,\k_p)_{\rho_k}^0$.

\begin{proposition}\label{prop:mingen}
    Let $m$ be such that the Hilbert series of $Q_m(D_n,\k_p)_{\rho_k}$ is different from the Hilbert series of $Q_m(D_n,\k_0)_{\rho_k}$, and let $a$ be the maximal possible value such that \eqref{eq:odddet}, \eqref{eq:evendet1}, or \eqref{eq:evendet2} holds. Then the Frobenius-type quasi-invariant constructed in Proposition \ref{prop:renxu} corresponding to $a$ is the minimal degree generator of $Q_m(D_n,\k_p)_{\rho_k}$.
\end{proposition}

\begin{proof}
    Let us first assume that $n$ is even. Then for any $m'=(m'_0,m'_1)\in\Z_{\geq -1}^2$, let $g$ be the Frobenius-type quasi-invariant constructed in $Q_m(D_n,\k_p)_{\rho_k}^0$ corresponding to $m'$. Let $m''$ be defined as in Proposition \ref{prop:renxu}, and let $f$ be the polynomial used to construct $g$ as in Lemmas \ref{lem:frob} and \ref{lem:twistedfrob}.
    
    Let
    \begin{align*}
        S_{m',a}=\Bigg\{m=(m_0,m_1)\in\Z_{\geq 0}^2\mid&\frac{m_0'+m_1'}{2}+\frac{\Fr^a(k)}{n}<\frac{\mu}{2p^a}<\frac{m_0'+m_1'}{2}+\frac{n-\Fr^a(k)}{n}-\frac{1}{p^a},\\
        &\frac{m_0'-m_1'}{2}-\frac{1}{2}+\frac{\Fr^a(k)}{n}+\frac{1}{2p^a}<\frac{\delta}{2p^a}<\frac{m_0'-m_1'}{2}+\frac{1}{2}-\frac{\Fr^a(k)}{n}-\frac{1}{2p^a}\Bigg\}
    \end{align*}
    if $m'_0+m'_1$ is even and
    \begin{align*}
        S_{m',a}=\Bigg\{m=(m_0,m_1)\in\Z_{\geq 0}^2\mid&\frac{m_0'+m_1'}{2}+\frac{1}{2}-\frac{\Fr^a(k)}{n}<\frac{\mu}{2p^a}<\frac{m_0'+m_1'}{2}+\frac{1}{2}+\frac{\Fr^a(k)}{n}-\frac{1}{p^a},\\
        &\frac{m_0'-m_1'}{2}-\frac{\Fr^a(k)}{n}+\frac{1}{2p^a}<\frac{\delta}{2p^a}<\frac{m_0'-m_1'}{2}+\frac{\Fr^a(k)}{n}-\frac{1}{2p^a}\Bigg\}
    \end{align*}
    if $m'_0+m'_1$ is odd. In the proof of Proposition \ref{prop:renxu}, we proved that $g$ has degree less than the minimal degree generator of $Q_m(D_n,\k_0)_{\rho_k}$ if $m\in S_{m',a}$. From the inequalities used in the proof, it is also clear that for integer-valued $m$, if $g$ has degree less than the minimal degree generator of $Q_m(D_n,\k_0)_{\rho_k}$, then $m\in S_{m',a}$.

    Let us now prove the proposition by induction on $\mu$. We only need to consider $m$ in some $S_{m',a}$, and we may assume $a$ is maximal. If $m$ is such that $b=0$ and at least one of $(m_0-1,m_1),(m_0,m_1-1)$ is in $S_{m',a}$, then without loss of generality, assume it is $(m_0-1,m_1)$. Then we have $b=0$ for $(m_0-1,m_1)$ as well, and by induction, the minimal degree generator of $Q_{(m_0-1,m_1)}(D_n,\k_p)_{\rho_k}$ is either $g$ or a Frobenius-type quasi-invariant $g'$ corresponding to some $a'>a$. Assume for contradiction that we have the latter. Then by Proposition \ref{prop:freeness}, we have $g=eg'$ for some $e\in\k_p[\h]^{D_n}$. Note that $g'w_0^2$ is a Frobenius-type quasi-invariant, and if it has degree less than $\frac{n\mu}{2}+k_\mu$ then it violates the maximality of $a$. But then we have $\deg g<\deg g'+\deg w_0^2=\deg g'+n$, so $\deg e<n$. Then up to constant factors, $e$ is a power of $r^2$. If $e$ is not constant, this implies $f$ is divisible by $r^2$.
    
    Assume for ease of notation that $m_1'\geq 0$. Let $k'=\Fr^a(k)$ if $m_0'\geq 0$ and $k'=\frac{n}{2}-\Fr^a(k)$ if $m_0'<0$, and define $m''=m'$ if $m'_0\geq 0$ and $m''=(0,m'_1)$ if $m'_0=-1$. If $f$ is divisible by $r^2$, then by Lemma \ref{lem:sympolys}, $f$ is not a generator of $Q_{m''}(D_n,\k_p)_{\rho_{k'}}$. In particular, by induction we have a polynomial $f'$ and some $a',b'$ such that either $f'^{p^{a'}}w_0^{2b'_0}w_1^{2b'_1}$, $f'^{p^{a'}}w_0^{2b'_0+1}w_1^{2b'_1}$, or $f'^{p^{a'}}w_0^{2b'_0}w_1^{2b'_1+1}$ is a lower degree polynomial in $Q_{m''}(D_n,\k_p)_{\rho_{k'}}$. In the first case, we obtain that $f'^{p^{a+a'}}w_0^{2p^ab'_0}w_1^{2p^ab'_1}$ is a Frobenius-type quasi-invariant in $Q_m(D_n,\k_p)_{\rho_k}$ of degree less than $g$, which violates the maximality of $a$. The second and third cases are similar. So we have a contradiction, and we have $e$ is constant, so $g=g'$ is the minimal degree generator of $Q_{(m_0-1,m_1)}(D_n,\k_p)_{\rho_k}$.
    
    We have $Q_m(D_n,\k_p)\subset Q_{(m_0-1,m_1)}(D_n,\k_p)$, so any generator of $Q_m(D_n,\k_p)_{\rho_k}$ has degree at least that of $g$. So $g$ is a minimal degree generator of $Q_m(D_n,\k_p)_{\rho_k}$, as desired.

    Now, suppose $b=0$ and both of $(m_0-1,m_1),(m_0,m_1-1)$ are not in $S_{m',a}$. If $m'_0+m'_1$ is even, then we must have
    \[\frac{\mu-1}{2p^a}=\frac{m'_0+m'_1}{2}+\frac{\Fr^a(k)}{n}-\epsilon\]
    for some $\epsilon<\frac{1}{2p^a}$. Then for any $a'<a$, we have
    \[\left\{\frac{\mu-1}{2p^{a'}}\right\}=\left\{p^{a-a'}\left(\frac{\Fr^a(k)}{n}-\epsilon\right)\right\}=\left\{\pm\frac{\Fr^{a'}(k)}{n}-p^{a-a'}\epsilon\right\}.\]
    Note that $p^{a-a'}\epsilon<\frac{1}{2p^{a'}}$, so one sees that $\mu-1$ is not in either of the ranges in \eqref{eq:evendet1} or \eqref{eq:evendet2}. In particular, $(m_0-1,m_1)$ does not satisfy either of the conditions. A similar argument yields the same conclusion when $m'_0+m'_1$ is odd.

    Then the minimal degree generator of $Q_{(m_0-1,m_1)}(D_n,\k_p)_{\rho_k}^0$ has degree $\frac{n(\mu-1)}{2}+\frac{n}{2}-k_\mu=\frac{n\mu}{2}-k_\mu$. Let $g'$ be the minimal degree generator of $Q_m(D_n,\k_p)_{\rho_k}^0$. Then we have $\deg g'\geq\frac{n\mu}{2}-k_\mu$, and by Proposition \ref{prop:freeness}, there is no other generator of $Q_m(D_n,\k_p)_{\rho_k}^0$ whose degree is at most that of $g$. So since $g\in Q_{(m_0-1,m_1)}(D_n,\k_p)_{\rho_k}^0$, we must have $g=eg'$ for some homogeneous $e\in\k_p[\h]^{D_n}$. But we have $\deg(e)\leq2k_\mu<n$, so up to constant factors, $e$ is a power of $r^{2}$. If $e$ is not constant, this would imply $g$ is divisible by $r^2$, which would imply $f$ is divisible by $r^2$. By the same argument as in the previous case, this cannot happen, so $g=g'$ is the minimal degree generator of $Q_m(D_n,\k_p)_{\rho_k}^0$.

    Finally, assume that $b\neq 0$, and without loss of generality say we have $b_0>0$. Let us first consider when $m'_0,m'_1\geq 0$. Let $C_0,C_1$ be the conjugacy classes in $\Sigma$ such that $s_0\in C_0,s_1\in C_1$. Then we will show that there is some $s\in C_0$ such that $(1-s)f^{p^a}$ is not divisible by $\alpha_s(x)^{p^a(2m'_0+1)+2}$. Indeed, if there was no such $s$, then $f$ would be $(m_0'+1,m_1')$-quasi-invariant. In particular, since $\deg f$ is lower than the minimal degree generator of $Q_{(m_0'+1,m_1')}(D_n,\k_0)_{\rho_{n/2-\Fr^a(k)}}$, this would imply by induction that $f$ can be written as a $\k_p[\h]^{D_n}$-multiple of a Frobenius-type quasi-invariant. But then by the same argument as above, this contradicts the maximality of $a$. The same argument also applies for $C_1$.

    By Proposition \ref{prop:freeness} and by induction, the minimal degree generator of $Q_m(D_n,\k_p)_{\rho_k}^0$ can be written as $ef^{p^a}$ for some $e\in\k_p[\h]^{D_n}$. There is some $s\in C_0$ such that $(1-s)f^{p^a}$ is not divisible by $\alpha_s(x)^{p^a(2m'_0+1)+2}$, so since $(1-s)(ef^{p^a})$ is divisible by $\alpha_s(x)^{2m_0+1}$, $e$ is divisible by $\alpha_s(x)^{2b_0}$. But then since $e$ is $D_n$-invariant, it is divisible by $w_0^{2b_0}$. The same argument applied to $C_1$ then implies that we must have $e=w_0^{2b_0}w_1^{2b_1}$, so $g$ is the minimal degree generator of $Q_m(D_n,\k_p)_{\rho_k}^0$.

    When $m'_0<0$, the same argument as above shows that there is some $s\in C_1$ such that $(1-s)(f^{p^a}w_0)$ is not divisible by $\alpha_s(x)^{p^a(2m'_1+1)+2}$ since $w_0$ is $s$-invariant. We will show that there is also some $s\in C_0$ such that $(1-s)(f^{p^a}w_0)$ is not divisible by $\alpha_s(x)^3$. Indeed, we have $(1-s)(f^{p^a}w_0)=w_0((1+s)f)^{p^a}$, so this would imply $(1+s)f$ is divisible by $\alpha_s(x)$. But then $f$ must be divisible by $\alpha_s(x)$. If this holds for all $s\in C_0$, then we have $f=w_0f'$ for some $f'$. But then we have $f^{p^a}w_0=f'^{p^a}w_0^{p^a+1}$, so $f^{p^a}w_0$ can be written as a Frobenius-type quasi-invariant in two different ways. But the sets $S_{m',a}$ for fixed $a$ are disjoint, so this cannot happen. Now, the same argument as above shows that $g$ is the minimal degree generator of $Q_m(D_n,\k_p)_{\rho_k}^0$.

    For odd $n$, a parallel argument holds, where instead of considering neighbors such as $(m_0-1,m_1)$ and $(m_0,m_1-1)$, we only need to consider $m-1$. We do not run into any issues transferring arguments about Frobenius-type quasi-invariants of the form $f^{p^a}w_0^{2m_0+1}w_1^{2m_1}$ since we never have $m_0=-1$ or $m_1=-1$, so these quasi-invariants do not occur.
    
\end{proof}

In particular, we obtain
\begin{corollary}\label{cor:degs}
    If $n$ is odd and $m$ satisfies \eqref{eq:odddet} for some $a$, then let $a$ be maximal such that \eqref{eq:odddet} is satisfied. Then let
    \[d=np^a\left\lfloor\frac{m}{p^a}\right\rfloor+p^a\Fr^a(k)+2nb\]
    where $b=\max\left(0,m-\frac{p^a(2\lfloor\frac{m}{p^a}\rfloor+1)-1}{2}\right)$. Then we have
    \[\H(Q_m(D_n,\k_p)_{\rho_k})=\frac{2t^d+2t^{2nm+n-d}}{(1-t^2)(1-t^n)}.\]
    
    If $n$ is even and $m$ satisfies \eqref{eq:evendet1} or \eqref{eq:evendet2} for some $a$, then let $a$ be maximal such that one of the inequalities is satisfied. Then let 
    \[d=\begin{cases}
        np^a\left\lfloor\frac{\mu}{2p^a}\right\rfloor+p^a\Fr^a(k)+n(b_0+b_1)&\text{\eqref{eq:evendet1} is satisfied}\\
        np^a\left\lfloor\frac{\mu}{2p^a}+\frac{1}{2}\right\rfloor-p^a\Fr^a(k)+n(b_0+b_1)&\text{\eqref{eq:evendet2} is satisfied}
    \end{cases}\]
    where $b$ is defined as in Proposition \ref{prop:renxu}. Then we have
    \[\H(Q_m(D_n,\k_p)_{\rho_k})=\frac{2t^d+2t^{n\mu+n-d}}{(1-t^2)(1-t^n)}.\]
\end{corollary}
\begin{proof}
    For the even case, combine Proposition \ref{prop:mingen} and the degree calculation from the proof of Proposition \ref{prop:renxu}. For the odd case, plug $m=\frac{\mu}{2}$ into the formula for \eqref{eq:evendet1}.
\end{proof}

\section*{AI Disclosure}
The results in Section \ref{sec:harmonic} were obtained using ChatGPT 5.4-5.5 Pro. The paper was checked by Claude Opus 5.5 and ChatGPT 6.1 Sol. Significant revisions and corrections by Claude Opus 5.5 to Proposition \ref{prop:freeness} and Proposition \ref{prop:renxu} and by ChatGPT 6.1 Sol to Lemma \ref{lem:odddet}, Lemma \ref{lem:evendet}, and Proposition \ref{prop:evenhyper} were used. Figure \ref{figure} was made using ChatGPT 6.1 Sol. Various typos were also caught by ChatGPT 6.1 Sol.

\section*{Acknowledgements}

The author would like to thank Pavel Etingof for his guidance and many useful discussions. This work is based upon work supported by the National Science Foundation Graduate Research Fellowship under Grant No.~2141064. Any opinions, findings, and conclusions or recommendations expressed in this material are those of the author(s) and do not necessarily reflect the views of the National Science Foundation.

\printbibliography

@article{Ren_2020,
   title={Quasi-Invariants in Characteristic p and Twisted Quasi-Invariants},
   journal={Symmetry, Integrability and Geometry: Methods and Applications},
   publisher={SIGMA (Symmetry, Integrability and Geometry: Methods and Application)},
   author={M. Ren and X. Xu},
   year={2020}
}

@article{felder2003action,
  title={Action of Coxeter groups on m-harmonic polynomials and Knizhnik-Zamolodchikov equations},
  author={G. Felder and A. Veselov},
  journal={Moscow Mathematical Journal},
  year={2003}
}

@article{feigin,
  title={Quasi-invariants of Coxeter groups and m-harmonic polynomials},
  author={M. Feigin and A. Veselov},
  journal={International Mathematics Research Notices},
  year={2002}
}

@article{etingof2002lectures,
  title={Lectures on quasi-invariants of Coxeter groups and the Cherednik algebra},
  author={P. Etingof and E. Strickland},
  journal={arXiv preprint math/0204104},
  year={2002}
}

@article{wang,
    author = {F. Wang},
    title = {Toward explicit Hilbert series of quasi-invariant polynomials in characteristic $p$ and $q$-deformed quasi-invariants},
    journal = {New York Journal of Mathematics},
    year = {2023},
    volume = {29}
}

@article{wangyee,
    author = {F. Wang and E. Yee},
    title = {Hilbert series of $S_3$-quasi-invariant polynomials in characteristics $2,3$},
    journal={Symmetry, Integrability and Geometry: Methods and Applications},
    year = {2025},
    volume = {21}
}

@misc{judy,
author = {S. Li},
title = {\textup{``Hilbert Series of Quasi-invariant Polynomials of $D_n$ in
Characteristic $p$"}},
howpublished = {Research Science Institute},
url = {https://math.mit.edu/documents/rsi/2025Li.pdf},
year = {2025}
}

@article{eg,
author = {P. Etingof and V. Ginzburg},
title = {On $m$-quasiinvariants of a Coxeter group},
journal = {Moscow Mathematical Journal},
year = {2002},
volume = {2}
}

@article{feigindihedral,
author = {M. Feigin},
title = {Quasi-invariants of dihedral systems},
journal = {Mathematical Notes},
year = {2004},
volume = {76}
}

@article{beg,
AUTHOR = {Yu. Berest and P. Etingof and V. Ginzburg},
TITLE = {{Cherednik algebras and differential operators on quasi-invariants}},
VOLUME = {118},
JOURNAL = {Duke Mathematical Journal},
PUBLISHER = {Duke University Press},
YEAR = {2003}
}

@misc{stacks-project,
  author       = {The {Stacks Project Authors}},
  title        = {Stacks Project},
  howpublished = {\url{https://stacks.math.columbia.edu}},
  year         = {2018},
}
\end{document}